\documentclass[aos,MSNbibl,seceqn,nameyear,dvips]{arximspdf}
\usepackage{graphicx}
\doi{10.1214/13-AOS1196} 
\volume{42}
\issue{2}
\pubyear{2014}
\firstpage{625}
\lastpage{653}

\makeatletter
\newcommand{\rrvert}{\vert}
\newcommand{\llvert}{\vert}
\newtheorem{teo}{Theorem}[section]
\newtheorem{lem}[teo]{Lemma}
\newtheorem{prop}[teo]{Proposition}
\newtheorem{cor}[teo]{Corollary}
\newproclaim{rem}[teo]{Remark}
\makeatother

\begin{document}
\begin{frontmatter}

\title{Convergence of linear functionals of the Grenander estimator
under misspecification}
\runtitle{Grenander under misspecification}

\begin{aug}
\author[A]{\fnms{Hanna} \snm{Jankowski}\corref{}\thanksref{t1}\ead[label=e1]{hkj@yorku.ca}}
\runauthor{H. Jankowski}
\affiliation{York University}
\address[A]{Department of Mathematics and Statistics\\
York University\\
4700 Keele Street\\
Toronto, Ontario M3J 1P3\\
Canada\\
\printead{e1}} 
\end{aug}
\thankstext{t1}{Supported in part by an NSERC Discovery Grant.}

\received{\smonth{8} \syear{2012}}
\revised{\smonth{12} \syear{2013}}

%
\begin{abstract}
Under the assumption that the true density is decreasing, it is well
known that the
Grenander estimator converges at rate $n^{1/3}$ if the true density is curved
[\textit{Sankhy\=a Ser. A} \textbf{31} (1969) 23--36]
and at rate $n^{1/2}$ if the density is flat
[\textit{Ann. Probab.} \textbf{11} (1983) 328--345;
\textit{Canad. J. Statist.} \textbf{27} (1999) 557--566].
In the case that the true density is misspecified, the results of
Patilea [\textit{Ann. Statist.} \textbf{29} (2001) 94--123]
tell us that the
global convergence rate is of order $n^{1/3}$ in Hellinger distance.
Here, we show that the local convergence rate is~$n^{1/2}$ at a point
where the density is misspecified. This is not in contradiction with the results
of Patilea [\textit{Ann. Statist.} \textbf{29} (2001) 94--123]:
the global convergence rate simply comes from locally curved \emph{well-specified} regions.
Furthermore, we study global convergence under misspecification by
considering linear functionals.
The rate of convergence is $n^{1/2}$ and we show that the limit is made
up of two independent terms:
a mean-zero Gaussian term and a second term (with nonzero mean) which
is present only if the
density has well-specified locally flat regions.
\end{abstract}

%

\end{frontmatter}

\section{Introduction}\label{sec1}

Shape-constrained nonparametric maximum likelihood estimators
provide an intriguing alternative to kernel-based density estimators.
For example, one can compare the standard histogram with the Grenander
estimator for a decreasing density. Rules exist to pick the bandwidth
(or bin width)
for the histogram to attain optimal convergence rates, cf. \citet
{AllNPstat}.
On the other hand, the Grenander estimator gives a piecewise constant
density, or
histogram, but the bin widths are now chosen completely automatically
by the estimator.
Furthermore, the bin widths selected by the Grenander estimator are
naturally locally
adaptive [\citet{birge}; \citet{cator}]. Similar
comparisons can also be made between the
log-concave nonparametric MLE and the kernel density estimator with,
say, the Gaussian kernel.

The Grenander estimator was first introduced in \citet{ulf2} and
has been considered
extensively in the literature since then. A recent review of the
history of the problem
appears in \citet{Linfinity}. The latter paper establishes that
the Grenander estimator
converges to a true strictly decreasing density at a rate of $(n/\log
n)^{1/3}$ in the
$L_\infty$ norm. Other rates have also been derived over the years,
most notably,
convergence at a point at a rate of $n^{1/3}$ if the true density is
locally strictly
decreasing [\citet{PrakasaRao}; \citet{Piet85}] and at a
rate of $n^{1/2}$ if the
true density is locally flat [\citet{piet83}; \citet
{CarolanDykstraCJS}].

As noted in \citet{CSS}; \citet{DSS}, the ``success story''
of maximum likelihood
estimators is their robustness. Namely, let $\mathcal{F}$ denote the
space of
decreasing densities on ${\mathbb R}_+$. Next, let $f_0$ denote the
true density
and $\hat f_0$ denote the density closest to $f_0$ in the
Kullback--Leibler sense. That is,
%
\begin{equation}
\label{linedefKL} \hat f_0 = \mathop{\operatorname{argmin}}_{g\in\mathcal{F}}
\int_0^\infty f_0(x) \log
\frac{f_0(x)}{g(x)} \,dx.
\end{equation}
We will call the density $\hat f_0$ the KL projection density of
$f_0$, or the
KL projection for short. Note that if $f_0\in\mathcal{F}$ then
$\hat f_0 = f_0$.
\citet{patilea} showed that the density $\hat f_0$ exists,
and that the Grenander estimator converges to $\hat f_0$ when the
observed samples come from the true density $f_0$, regardless if
$f_0\in\mathcal{F}$.
Similar results were proved for the log-concave maximum likelihood
estimator in
\citet{CuleSam}; \citet{CSS}; \citet{DSS};
\citet{BJR}.

\begin{figure}

\includegraphics{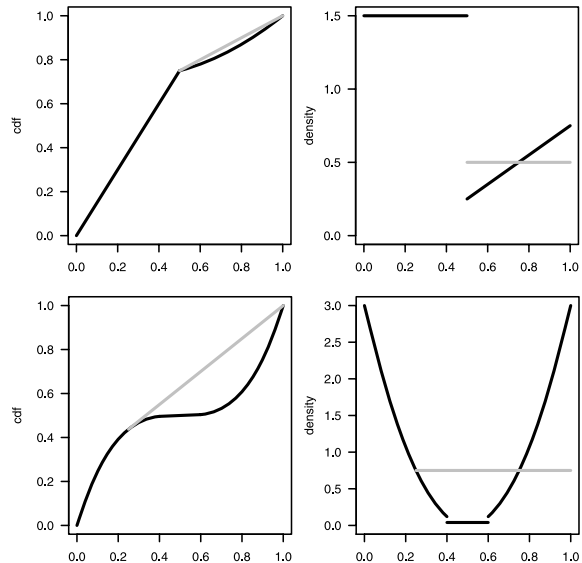}

\caption{Two examples of $f_0$ and $\hat f_0 = \operatorname
{gren}(F_0)$. The
two left panels show the c.d.f. and density, for example, (\protect\ref{lineeg1})
while the two right panels\vspace*{1pt} show the c.d.f. and density, for example,
(\protect\ref{lineeg2}).
$F_0$ (resp.,~$f_0$) is shown in black, and $\widehat F_0$ (resp.,~$\hat f_0$)
is shown in gray, but only if different from the truth (namely $F_0$
and $f_0$, resp.).} \label{figmisseg}
\end{figure}

In order to understand the local behavior of the Grenander estimator when
\mbox{$f_0\notin\mathcal{F}$}, we first need to define regions where $f_0$
is considered to be
miss- and well-specified. Let $\widehat F_0$ denote the cumulative
distribution function
of $\hat f_0$ defined in~(\ref{linedefKL}). The regions where
$\widehat F_0 \neq F_0$
are then the regions where $f_0$ is misspecified, and $f_0$ is
considered to be
well specified otherwise. Note that, if $f_0$ is misspecified in a
region, it may still
be decreasing on some portion of this region; see, for example,
Figure~\ref{figmisseg}.

%

Let $\hat f_n$ denote the Grenander estimator of a decreasing density.
We show here that at a point where the density is misspecified the rate of
convergence of $\hat f_n$ to $\hat f_0$ is $n^{1/2}$, and we also
identify the limiting distribution. This is not in contradiction with
the results of
\citet{patilea}: the slower $n^{1/3}$ global convergence rate
simply comes
from locally curved {well-specified} regions. To be more specific, if
the density
$f_0$ is misspecified at a point, then $\widehat F_0$ must be linear (and
$\hat f_0$ is flat), and in regions where $\hat f_0$ is flat
the rate of
convergence is $n^{1/2}$. In fact, the $n^{1/2}$ rate holds at all flat regions
of $\hat f_0$, irrespective of whether these are miss- or well-specified.
The complete result is given in Section~\ref{secpointwise}, where
some properties
of $\hat f_0$ are also discussed. 

Next, we consider convergence of linear functionals. Let
%
\begin{equation}
\label{defmug} \hat\mu_0(g) = \int_0^\infty
g(x) \hat f_0(x)\,dx\quad\mbox{and}\quad \hat\mu_n(g) =
\int_0^\infty g(x) \hat f_n(x)\,dx.
\end{equation}
In Section~\ref{seclinear}, we show that $n^{1/2}(\hat\mu
_n(g)-\hat\mu_0(g)) = O_p(1)$,
and we again identify the limiting distribution. Notably, the limit is
made up of two
\emph{independent} terms: a~mean-zero Gaussian term and a second term
with nonzero mean. Furthermore, the second term is present only if the
density has {well-specified locally flat} regions. Our results apply to
a wide
range of KL projections with both strictly curved and flat regions.
The work in the strictly curved case follows from the rates of
convergence of
$\widehat F_n(y)=\int_0^y \hat f_n(y)\,dy$ to the empirical distribution
function established in \citet{KW1}. However, as mentioned above,
this is only for
the well-specified regions of~$f_0$. A related work here is that of
\citet{KulikovLopuhaaJTP}, who consider functionals in the
strictly curved case
but at the distribution function level.

In Section~\ref{secentropy}, we go beyond the linear setting, and consider
convergence of the entropy functional in the misspecified case. The
limit in this
case is Gaussian, irrespective of the local properties of $\hat f_0$.
Most proofs appear in Section~\ref{secproofs} and some technical
details are given in \citet{suppmat}.\vadjust{\goodbreak} Throughout, our results are
illustrated by reproducible simulations.
Code for these is available online at \href{http://www.math.yorku.ca/~hkj/}{www.math.yorku.ca/~hkj/}.

To our best knowledge, previous work on rates of convergence under
misspecification in the shape-constrained context is limited to the
rates established
in \citet{geerbook} and \citet{patilea}, as well as the
more recent results of \citet{BJR}.
In \citet{BJR}, the pointwise asymptotic distribution under
misspecification was derived
for the log-concave probability mass function.

The implications of the new results obtained here are as follows.
First, we now
understand that $\hat f_0$ will be made up of local well specified and
misspecified\vspace*{1pt} regions, and that the rate of convergence in the
misspecified regions
is always $n^{1/2}$. We conjecture that this type of behavior will be
seen in other situations,
such as the log-concave\vspace*{1pt} setting for $d=1$. That is, the rate of
convergence in misspecified
regions will be $n^{1/2}$ whereas in well-specified regions the rate
of convergence will
depend on whether locally the density lies on the boundary or the
interior of the underlying space.
In the log-concave $d=1$ case, this ``interior'' rate is known to be
$n^{2/5}$ [\citet{BRW}].
The interesting case of $d>1$ is more mysterious though, as the
relationship between
the slower boundary points and faster interior points is harder to identify.

Secondly, we show that linear functionals (as well as the nonlinear
entropy functional)
converge at rate $n^{1/2}$, and we also conjecture that this behavior
will continue to
hold for other shape constraints. Let $\mu_0(g)=\int_0^\infty
g(x)f_0(x)\,dx$. Our results show that
%
\begin{equation}
\label{linepower} \sqrt{n} \bigl(\hat\mu_n(g)-\mu_0(g)
\bigr) = O_p(1) + \sqrt{n} \bigl(\hat\mu_0(g)-
\mu_0(g) \bigr).
\end{equation}
%
Therefore,
global rates of \emph{divergence} are $n^{1/2}$ for linear functionals
in the misspecified case. A similar statement also holds for the
entropy functional,
and here the random $O_p(1)$ term is always Gaussian. Such results are
well understood in parametric settings, and are key in power calculations.
The exact conditions necessary for (\ref{linepower}) to hold are
given in
Section~\ref{seclinear} for $\mu_0(g)$ and in Section~\ref{secentropy}
for the entropy. Our work can also be easily extended to locally misspecified
settings such as those studied in \citet{LeCam}.

Lastly, the fact that the limiting distribution of the linear
functional $\hat\mu_n(g)$ 
depends on properties of $\hat f_0$, whereas the limiting
distribution of the
entropy functional is always Gaussian, makes the entropy functional potentially
more appealing in terms of testing procedures. Hypothesis testing based
on functionals was considered, for example, in \citet{CSS} and
\citet{trace}. The latter reference develops the ``trace test''
which depends on a nearly linear functional, the variance.
Both, however, are developed in the context of log-concavity, and it
would be
of great interest to extend the results presented here to that setting,
particularly for higher dimensions.


\section{The Kullback--Leibler projection and pointwise convergence under misspecification}\label{secpointwise}

Properties of the KL projection onto the space of log-concave densities
were studied in \citet{DSS}. When \mbox{projecting} onto the space of decreasing
densities, the behavior is a little easier to characterize.

%
\begin{teo}[{[\citeauthor{patileaphd} (\citeyear{patileaphd,patilea})]}]\label{teopatilea}
Let\vspace*{1pt} $f_0$ be a density with support on~$[0,\infty)$ with
$F_0(x)=\int_0^x f_0(u)\,du$. Let $\widehat F_0$ denote the least
concave majorant of $F_0$. Then the left derivative of $\widehat F_0,
\hat f_0$, satisfies the inequality \mbox{$\int\log\frac{\hat f_0}{f} \,dF_0
\geq0$}, for all decreasing densities $f$.\vadjust{\goodbreak}
\end{teo}

%
\begin{rem}
The density $\hat f_0$ satisfying $\int\log\frac{\hat f_0}{f} \,dF_0
\geq0$ for all $f\in\mathcal{F}$ is called the
``pseudo-true'' density by \citet{patilea}. If we additionally
assume that $\sup_{f\in\mathcal{F}} \int\log f \,dF_0$ and $\int
\log f_0 \,dF_0$ are both finite, then this $\hat f_0$ is also the
unique minimizer of the Kullback--Leibler divergence
\[
\hat f_0 = \mathop{\operatorname{argmin}}_{f \in\mathcal{F}} \int\log
\frac{f_0}{f} \,dF_0.
\]
See \citet{patilea}, page~95, for more details. In what follows,
we continue to refer to $\hat f_0$ as defined in Theorem~\ref{teopatilea}, as the KL projection, even if it comes from the slightly
more general definition of \citet{patilea}.
\end{rem}

Thus, in our setting, we have a complete graphical representation
of the distribution function $\widehat F_0$ of the KL projection.
This representation makes it possible to calculate $\hat f_0$ in
many cases.\vspace*{1pt}
It also allows us to easily visualize the various $F_0$ which yield the
same $\hat f_0$.
Moreover, the representation is key in understanding the behavior of
the estimator,
both on the finite sample and asymptotic levels. Therefore, for a
function $g$
we define the operator $\operatorname{gren}(g)$ to denote the (left)
derivative of the least
concave majorant of $g$. When the least concave majorant is restricted
to a set
$[a,b]$, we will write $\operatorname{gren}_{[a,b]}(g)$.

Let $\mathcal{S}_0$ denote the support of $f_0$. We write
$\mathcal{S}_0 = \mathcal{M} \cup\mathcal{W}$, where $\mathcal
{M}=\{x \geq0\dvtx  \widehat F_0(x) > F_0(x)\}$
and $\mathcal{W}=\{x \geq0\dvtx  \widehat F_0(x) = F_0(x)\}$.
Since\vspace*{1pt} $f_0$ is a density, it follows that $F_0$ is continuous, as is
$\widehat F_0$ and, therefore, $\mathcal{W}$ is a closed set and
$\mathcal{M}$ is open.
For a fixed point $x_0\in\mathcal{M}$, we thus know that $x_0$ lies
in\vspace*{1pt} some open interval.
Indeed, let $a_0=\sup\{x<x_0\dvtx  \widehat F_0(x)=F_0(x)\}$ and
$b_0=\inf\{x>x_0\dvtx  \widehat F_0(x)=F_0(x)\}$. Then $x_0\in(a_0,b_0)$
with $(a_0,b_0)\subset\mathcal{M}$.

Two examples are given in Figure~\ref{figmisseg}. For the first
example, we have
%
\begin{eqnarray}
\label{lineeg1} f_0(x)&=& \cases{1.5, &\quad$x \in[0, 0.5]$,
\vspace*{2pt}\cr
x-0.25, &\quad$x \in(0.5,1]$.}
\end{eqnarray}
Here, $\mathcal{M}=(0.5,1)$ and $\mathcal{W}=[0, 0.5]\cup\{1\}$. For
the second example, we have
%
\begin{eqnarray}
\label{lineeg2} f_0(x)&=& \cases{ 12(x-0.5)^2, &\quad$x
\in[0,0.4]\cup[0.6,1]$,
\vspace*{2pt}\cr
0.04, &\quad$x \in(0.4,0.6)$.}
\end{eqnarray}
Here, $\mathcal{M}=(0.25,1)$ and $\mathcal{W}=[0, 0.25]\cup\{1\}$.

The next proposition gives some additional properties of the KL projection.

%
\begin{prop}\label{propKLprop}
The density, $\hat f_0$, satisfies the following:
\begin{longlist}[(2)]
\item[(1)] Fix $x_0\in\mathcal{M}$ and define $a_0, b_0$ as above.
Then $b_0<\infty$, and $\hat f_0$ is constant on $(a_0,b_0]$
and satisfies the mean-value property
\[
\hat f_0(x_0) = \frac{1}{b_0-a_0} \int
_{a_0}^{b_0} f_0(x)\,dx.
\]
\item[(2)] Suppose that $\int_0^\infty f_0^2(x)\,dx<\infty$.
Then
\[
\hat f_0 = \mathop{\operatorname{argmin}}_{g\in\mathcal{F}} \int_0^\infty
\bigl(g(x)-f_0(x)\bigr)^2 \,dx.
\]

\item[(3)] For any increasing function $h(x)$,
$\int_0^\infty h(x) \hat f_0(x)\,dx \leq\int_0^\infty h(x) f_0(x)\,dx$.
%
\item[(4)] Let $g_0\in\mathcal{F}$ and let $G_0(y)=\int_0^y g_0(x)\,dx$. Then
\begin{eqnarray*}
\sup_{x\geq0} \bigl|\widehat F_0(x) -
G_0(x)\bigr| &\leq& \sup_{x\geq0} \bigl|F_0(x) -
G_0(x)\bigr|.
\end{eqnarray*}
\end{longlist}
\end{prop}

Point (3) above tells us that if $g$ is increasing then
$\mu_0(g)\geq\hat\mu_0(g)$. Point (4) is Marshall's lemma
[\citet{marshall}].
The proof of Proposition \ref{propKLprop} appears in \citet{suppmat}.

Suppose that $X_1, \ldots, X_n$ are independent
and identically distributed with density $f_0$ on ${\mathbb R}_+=[0,
\infty)$.
Let $\hat f_n$ denote the Grenander estimator of a decreasing density
\[
\hat f_n = \mathop{\operatorname{argmax}}_{g \in\mathcal{F}} \int\log
g(x) \,d{ \mathbb F}_n(x),
\]
where $\mathcal{F}$ denotes the class of decreasing densities on
${\mathbb R}_+$, and ${\mathbb F}_n(x) = n^{-1} \sum_{i=1}^n
1_{(-\infty, x]}(X_i)$
denotes the empirical distribution function. The next theorem is our
first main result.

%
\begin{teo}\label{teomisslocal}
Fix\vspace*{1.5pt} a point $x_0\in\mathcal{M}$, and let $[a,b]$ denote the
largest interval~$I$ containing $x_0$ such that $\widehat F_0(x)$ is
linear on $I$.
Let ${\mathbb U}$ denote a standard Brownian bridge process on $[0,1]$,
and let ${\mathbb U}_{F_0}(x)={\mathbb U}(F_0(x))$ for $x\in\mathcal
{S}_0$. Then
\begin{eqnarray*}
\sqrt{n} \bigl(\hat{f}_n(x_0)-\hat f_0(x_0)
\bigr) &\Rightarrow& \operatorname{gren}_{[a,b]} \bigl({\mathbb
U}^{\mathrm{mod}}_{F_0} \bigr) (x_0),
\end{eqnarray*}
where
\begin{eqnarray*}
{\mathbb U}^{\mathrm{mod}}_{F_0}(u) &=& \cases{ {\mathbb
U}_{F_0}(u), &\quad$u \in[a,b] \cap\mathcal{W}$,
\vspace*{2pt}\cr
-\infty, &\quad$u
\in[a,b] \cap\mathcal{M}$.}
\end{eqnarray*}
If it happens that $[a,b]\cap\mathcal{W} = \{a,b\}$, then
\begin{eqnarray*}
\sqrt{n} \bigl(\hat{f}_n(x_0)-\hat f_0(x_0)
\bigr) &\Rightarrow& \sigma Z,
\end{eqnarray*}
where $Z$ is a standard normal random variable and
\begin{eqnarray*}
\sigma^2 &=& \hat f_0(x_0) \biggl[
\frac{1}{b-a} -\hat f_0(x_0) \biggr].
\end{eqnarray*}
\end{teo}

Recall\vspace*{1pt} that \citet{patilea}, Corollary~5.6, shows that the rate of
convergence (in Hellinger distance) of $\hat{f}_n$ to $\hat f_0$ is $n^{1/3}$.
The above theorem shows that the \emph{local} rate of convergence
will be $\sqrt{n}$ where the KL projection is flat. When the KL density
is curved, the KL density and true density are actually equal, and
hence the convergence rate from the correctly specified case applies.
The next formulation of the limiting process is similar to that of
\citet{CarolanDykstraCJS} for a density with a flat region on $[a,b]$.

%
\begin{rem}\label{remmisslocal}
Let $p_0 =F_0(b)-F_0(a)= \widehat F_0(b)-\widehat F_0(a)$.
Since $\widehat F_0$ is linear on $[a,b]$ the limiting distribution may
also be expressed as
\begin{eqnarray*}
\operatorname{gren}_{[a,b]} \bigl({\mathbb U}^{\mathrm{mod}}_{F_0}
\bigr) (x_0) &=& \frac{1}{b-a} \biggl\{ Z+ \sqrt{p_0}
\operatorname{gren} \bigl({\mathbb U}^{\mathrm{mod}} \bigr) \biggl(
\frac{x_0-a}{b-a} \biggr) \biggr\},
\end{eqnarray*}
where $Z$ is a mean zero normal random variable with variance
$p_0(1-p_0)$, ${\mathbb U}$ is an independent standard Brownian bridge, and
\begin{eqnarray*}
{\mathbb U}^{\mathrm{mod}}(u) &=& \cases{ {\mathbb U}(u), &\quad$u \in
\bigl([a,b]\cap\mathcal{W}-a \bigr)/(b-a)$,
\vspace*{2pt}\cr
-\infty, &\quad$u \in \bigl([a,b]
\cap\mathcal{M}-a \bigr)/(b-a)$.}
\end{eqnarray*}
Notably, if $[a,b]\cap\mathcal{W} = \{a,b\}$, then $\operatorname
{gren}({\mathbb U}^{\mathrm{mod}})(u)=0$.
\end{rem}

%
\begin{figure}

\includegraphics{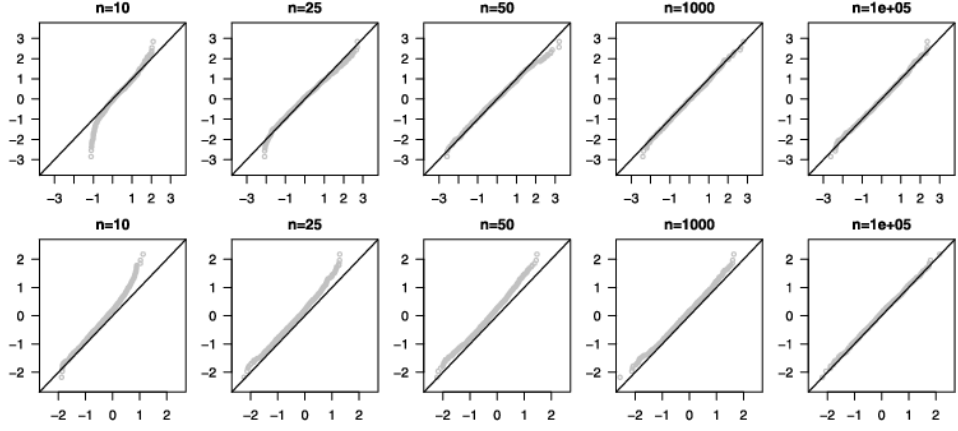}

\caption{Empirical quantiles of $\sqrt{n}(\hat f_n(x_0)-\hat
f_0(x_0))$ vs. the true quantiles of the limiting $N(0,\sigma^2)$
distributions at the point $x_0=0.75$ for $f_0$ given by (\protect\ref{lineeg1})
in the top row ($\sigma^2=3/4)$ and (\protect\ref{lineeg2}) in
the bottom row ($\sigma^2=7/16$). The sample size varies from $n=10$
to $n={}$100,000. The straight line goes through the origin and has
slope one. Each plot is based on $B=1000$ samples.}
\label{figmisspoint}
\end{figure}

Figure~\ref{figmisspoint} illustrates the theory. The convergence is
surprisingly fast, although it appears to be a little slower in the
second example (\ref{lineeg2}). We conjecture that this difference is
caused by the presence/absence of the strictly curved region of~$f_0$.


\begin{pf*}{Proof of Theorem \ref{teomisslocal}}
By the switching relation [\citet{BJPSW10}], we have
\begin{eqnarray*}
&& P \bigl(\sqrt{n} \bigl(\hat f_n(x_0)- \hat
f_0(x_0) \bigr) < t \bigr)
\\[-1pt]
&&\qquad = P \Bigl(\mathop{\operatorname{argmax}}_{z\geq0} \bigl\{{\mathbb
F}_n(z)- \bigl(\hat f_0(x_0)+n^{-1/2}t
\bigr)z \bigr\} < x \Bigr)
\\[-1pt]
&&\qquad = P \Bigl(\mathop{\operatorname{argmax}}_{z\geq0} \bigl\{\sqrt{n}
\bigl({\mathbb F}_n(z)-{\mathbb F}_n(a)-
\bigl(F_0(z)- F_0(a) \bigr) \bigr)
\\[-1pt]
&&\hspace*{83pt}{} + \sqrt{n}
\bigl(F_0(z)-F_0(a)- \hat f_0(x_0)
(z-a) \bigr)-tz \bigr\} < x \Bigr).
\end{eqnarray*}
We now look more closely at the ``second'' term. That is,
\begin{eqnarray*}
&& F_0(z)-F_0(a)-\hat f_0(x_0)
(z-a)
\\
&&\qquad = - \bigl\{\widehat F_0(z)-F_0(z) \bigr\} +
\bigl\{ \widehat F_0(z)-\widehat F_0(a)-\hat
f_0(x_0) (z-a) \bigr\},
\end{eqnarray*}
noting that $\widehat F_0(a)=F_0(a)$, since $a \in[a,b]\cap\mathcal{W}$.
On the other hand, for all $z \in[a,b] \cap\mathcal{M}$, we have
$\widehat F_0(z) > F_0(z)$. Furthermore, $\widehat F_0$
is concave with derivative $\hat f_0(x_0)$ [at any point $z\in
(a,b)$], and hence
\[
\widehat F_0(z)- \widehat F_0(a) - \hat
f_0(x_0) (z-a) \leq0
\]
for all $z \geq0$. For $z\in[a,b]\cap\mathcal{W}$ this is an
equality, and a strict inequality otherwise. Therefore, the weak limit of
\[
\sqrt{n} \bigl\{{\mathbb F}_n(z)-{\mathbb F}_n(a)-
\bigl(F_0(z)- F_0(a) \bigr) \bigr\} - \sqrt{n}
\bigl(F_0(z)-F_0(a)-\hat f_0(x_0)
(z-a) \bigr)
\]
is
${\mathbb U}_{F_0}^{\mathrm{mod}}(z)- {\mathbb U}_{F_0}^{\mathrm{mod}}(a) = {\mathbb
U}_{F_0}^{\mathrm{mod}}(z)- {\mathbb U}_{F_0}(a)$,
for all $z\in[a,b]$. For $z \notin[a,b]\cap\mathcal{W}$, the limit
of this
process is always $-\infty$ and, therefore, the maximum must occur
inside of $[a,b]$. By the argmax continuous mapping theorem
[\citet{jawbook}, Theorem 3.2.2, page 287],
\begin{eqnarray*}
&& P \bigl(\sqrt{n} \bigl(\hat f_n(x_0)-\hat
f_0(x_0) \bigr) < t \bigr)
\\
&&\qquad \rightarrow P \Bigl(
\mathop{\operatorname{argmax}}_{z\in[a,b]} \bigl\{ {\mathbb
U}_{F_0}^{\mathrm{mod}}(z) - tz \bigr\} < x \Bigr)
= P \bigl(\operatorname{gren}_{[a,b]} \bigl({\mathbb
U}_{F_0}^{\mathrm{mod}} \bigr) (x_0) < t \bigr)
\end{eqnarray*}
by switching again. When $[a,b]\cap\mathcal{W} = \{a,b\}$, then the
least concave majorant is simply the line joining ${\mathbb
U}_{F_0}(a)$ and ${\mathbb U}_{F_0}(b)$, with slope equal to
\[
\frac{{\mathbb U}_{F_0}(b)-{\mathbb U}_{F_0}(a)}{b-a},
\]
a Gaussian random variable with mean zero and variance
\begin{eqnarray*}
&& \frac{1}{(b-a)^2} \bigl(F_0(b)-F_0(a) \bigr) \bigl[1-
\bigl(F_0(b)-F_0(a) \bigr) \bigr]
\\
&&\qquad  =  \hat f_0(x_0) \biggl[
\frac{1}{b-a} -\hat f_0(x_0) \biggr].
\end{eqnarray*}\upqed
%
\end{pf*}

\begin{pf*}{Proof of Remark \ref{remmisslocal}}
Recall that $\widehat F_0$ is linear on $[a,b]$. Therefore, for $x\in[a,b]$,
we can write ${\mathbb U}(\widehat F_0(x))-{\mathbb U}(\widehat F_0(a))
= \frac{x-a}{b-a} {\mathbb W}+ {\mathbb V}(x)$,
where
\begin{eqnarray*}
{\mathbb W}&=& {\mathbb U} \bigl(\widehat F_0(b) \bigr)-{\mathbb U}
\bigl(\widehat F_0(a) \bigr),
\\
{\mathbb V}(x)&=& {\mathbb U} \bigl(\widehat F_0(x) \bigr)-{\mathbb
U} \bigl(\widehat F_0(a) \bigr) - \frac{\widehat F_0(x)-\widehat
F_0(a)}{\widehat F_0(b)-\widehat
F_0(a)} {\mathbb W}
\\
&=& {\mathbb U} \bigl(\widehat F_0(x) \bigr)-{\mathbb U} \bigl(
\widehat F_0(a) \bigr) - \frac
{x-a}{b-a} {\mathbb W}.
\end{eqnarray*}
Since all variables are jointly Gaussian, a careful calculation of the
covariances reveals that ${\mathbb W}$ and ${\mathbb V}(x)$ are
independent (also as processes),
and ${\mathbb W}$ is mean-zero Gaussian with variance $p_0(1-p_0)$. Furthermore,
\[
{\mathbb V}(s) \stackrel{d} {=} \sqrt{p_0} {\mathbb U} \biggl(
\frac
{s-a}{b-a} \biggr).
\]
This decomposition is similar to that of
\citet{ShorackWellnerBook}, Exercise~2.2.11, page~32.
Now, note that the Grenander operator satisfies
$\operatorname{gren}_{[a,b]}(g)(x) = \beta+ \frac{\gamma}{b-a}
\operatorname{gren}_{[0,1]}(h) (\frac{t-a}{b-a} )$ if\vspace*{1pt}
$g(t)=\alpha+\beta t + \gamma h (\frac{t-a}{b-a} )$.
It follows that
\begin{eqnarray*}
\operatorname{gren}_{[a,b]} ({\mathbb U}_{\widehat F_0} )
(x_0) &=& \frac{1}{b-a} Z+ \frac{\sqrt{p_0}}{b-a}\operatorname
{gren}({\mathbb U}) \biggl(\frac{x_0-a}{b-a} \biggr)
\end{eqnarray*}
with $Z, {\mathbb U}$ defined as in the remark.
The full result follows since,
${\mathbb U}^{\mathrm{mod}}_{F_0}(x)={\mathbb U}^{\mathrm{mod}}_{\widehat F_0}(x)=\frac
{x-a}{b-a} {\mathbb W}+ {\mathbb V}^{\mathrm{mod}}(x)$.
\end{pf*}

\section{\texorpdfstring{$\sqrt{n}$}{(Square root of n)}-convergence of linear functionals}\label{seclinear}

Consider a density $f_0$ with support $\mathcal{S}_0$ and let
$\hat f_0$ denote its KL projection.
We write $\mathcal{S}_0 = \mathcal{S}_c \cup\mathcal{S}_f$, where
$\mathcal{S}_c$ denotes
the portion of the support where $\hat f_0$ is curved and $\mathcal{S}_f$
denotes the portion of the support where $\hat f_0$ is flat. By
definition of $S_f$ as well as Proposition~\ref{propKLprop}, the KL
projection can be written as
%
\begin{eqnarray}
\label{defKLform} \hat f_0(x) &=& \sum_{j=1}^J
\hat q_j 1_{I_j} (x)
\end{eqnarray}
on $S_f$, where the intervals are disjoint and each is of the form
$I_j=(a_j, b_j]$.
Our results for linear functionals hold under the following assumptions:

\begin{longlist}[(C)]
\item[(S)] The support, $\mathcal{S}_0$, of $f_0$ is bounded.
\item[(C)] When\vspace*{1pt} the KL projection is curved, $\sup_{x \in\mathcal
{S}_c} |f'_0(x)| < +\infty$.
\item[(P)] The true density is strictly positive: $\inf_{x \in S_0}
f_0(x) >0$.
%
\item[(F)] 
When the KL projection is flat, $J $ is finite in (\ref{defKLform}).
\end{longlist}
Let $g\dvtx \mathcal{S}_0 \mapsto{\mathbb R}$ and define $\hat\mu_n(g)$
by (\ref{defmug}). Then we require that $g$ satisfy the following conditions:
\begin{longlist}[(G2)]
\item[(G1)] $\int_{\mathcal{S}_c} |g'(x)|\,dx <\infty$.\vspace*{1pt}
\item[(G2)] $g\in L_\beta(\mathcal{S}_f)$ for some $\beta>2$.
\end{longlist}

In order to state our main result for linear functionals, we need to
define the following functions: 
%
\begin{eqnarray}
g_j(u)&=& g \bigl((b_j -
a_j)u+ a_j \bigr),\qquad u \in[0,1],
\nonumber\\[-9pt]\label{linedefg}  \\[-9pt]
\bar{g}_j&=& (b_j-a_j)^{-1}
\int_{a_j}^{b_j} g(x)\,dx\nonumber
\end{eqnarray}
and
\begin{equation}\label{linedefg2}
\bar{g}(x)  =  \cases{ g(x), &\quad$x \in\mathcal{S}_c$,
\vspace*{2pt}\cr
\bar g_j, &\quad$x \in I_j, j=1, \ldots, J$.}
\end{equation}
Thus, $\bar{g}_1, \ldots, \bar{g}_J$ are the local
averages of the
function $g$, and each $g_j(u)$ is a localized version of $g$.

%
\begin{teo}\label{teolinear}
Suppose that the density $f_0$ satisfies conditions \textup{(S)}, \textup{(C)}, \textup{(P)}~and~\textup{(F)}.
Consider a function $g\dvtx \mathcal{S}_0 \mapsto{\mathbb R}$ which
satisfies conditions \textup{(G1)} and~\textup{(G2)}.
Let ${\mathbb U}, {\mathbb U}_1, \ldots, {\mathbb U}_J$ denote
independent Brownian bridges,
${\mathbb U}_{F_0}(x)={\mathbb U}(F_0(x))$, and define ${\mathbb
U}_j^{\mathrm{mod}}$ as in Theorem \ref{teomisslocal}.
Then
\begin{eqnarray*}
\sqrt{n} \bigl(\hat\mu_n(g)-\hat\mu_0(g) \bigr) &
\Rightarrow& \int_{\mathcal{S}_0} \bar{g}(x) \,d{\mathbb
U}_{F_0} (x)
\\
&&{} + \sum_{j=1}^J \sqrt{p_j}
\int_0^1 g_j(u)
\operatorname{gren} \bigl({\mathbb U}_{j}^{\mathrm{mod}} \bigr) (u)
\,du,
\end{eqnarray*}
where $p_j = F_0(b_j)-F_0(a_j) = \widehat F_0(b_j)-\widehat F_0(a_j)$.
Furthermore,
\[
\int_{\mathcal{S}_0} \bar{g}(x) \,d{\mathbb U}_{F_0} (x)
=\int_{\mathcal{S}_0} \bar{g}(x) \,d{\mathbb U}_{\widehat F_0} (x).
\]
Also, if $I_j\cap\mathcal{W} = \{a_j, b_j\}$, then $\operatorname
{gren}({\mathbb U}_{j}^{\mathrm{mod}})\equiv0$.
\end{teo}

%
\begin{figure}

\includegraphics{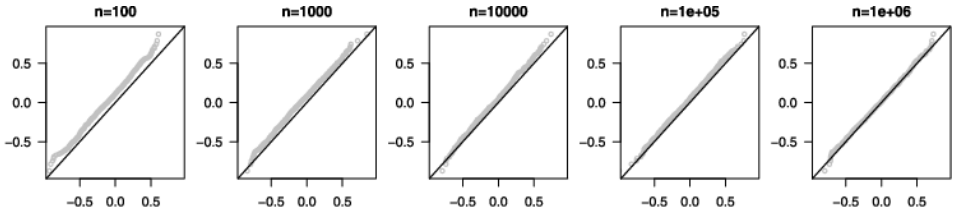}

\caption{Empirical quantiles of
$\sqrt{n}(\hat\mu_n(g)-\hat\mu_0(g))$ with $g(x)=x$ vs.
the true quantiles of the limiting $N(0,\sigma^2)$ distribution for
$f_0$ given in (\protect\ref{lineeg2}),
with $\sigma^2 \approx0.07032$.}\label{figmissmean}
\end{figure}

It follows that $\sqrt{n}(\hat\mu_n(g)-\hat\mu_0(g))$ will
converge to a Gaussian limit for true density (\ref{lineeg2}) but not
for (\ref{lineeg1}), as the latter has well-specified flat regions. A
simulation for (\ref{lineeg2}) is shown in Figure~\ref{figmissmean}.
The proof of Theorem~\ref{teolinear} is given in
Section~\ref{secproofs}. The simulations show that there appears a
systematic bias prior to convergence (the empirical quantiles appear on
the $x$-axis in Figure~\ref{figmissmean}, the negative bias
translates to a left-shift in the plot). The proof of Proposition~\ref{curved} shows that one source of the bias is the term $\sqrt{n}\int x
d(\widehat F_n - {\mathbb F}_n) \approx- \sqrt{n}\int(\widehat F_n -
{\mathbb F}_n) \leq0$. When $\mathcal{S}_0=\mathcal{S}_c$, this term
is the only source of bias, and from \citet{KW1}, it converges to
zero at a rate of at least $n^{1/6} (\log n)^{-2/3}$. Since (3) of
Proposition \ref{propKLprop} also holds at the empirical level,
similar behavior will be seen for all increasing functions $g$.

The results of Theorem \ref{teolinear} also show that $\sqrt{n}(\hat\mu
_n(g)-\hat\mu_0(g))$ is asymptotically normal
with variance $\operatorname{var}_{f_0}(\bar g(X))=\operatorname
{var}_{\hat f_0}(\bar g(X))$ if $\mathcal{S}_0$ has no \emph{well-specified flat} regions. Additionally, if $\mathcal{S}_0 =
\mathcal{S}_c$, then $\bar g(x)=g(x)$ and the model is well
specified. In this case, $\hat\mu_n(g)$ has the same asymptotic
distribution as the empirical estimator $n^{-1} \sum_{i=1}^n g(X_i)$
(see also Proposition~\ref{curved}). This shows that the maximum
likelihood estimator is asymptotically efficient, as in the strictly
curved case the family of decreasing densities is complete, and hence
the ``naive estimator'' $n^{-1} \sum_{i=1}^n g(X_i)$ is asymptotically
efficient [\citet{saramix}, Example 4.7].

Finally, we make a few comments on the assumptions required for Theorem~\ref{teolinear} to hold. 
The assumptions which we use on $S_c$ are \textup{(S)}, \textup{(P)} and \textup{(C)}. These are
quite standard assumptions in the literature for the strictly curved
setting; see, for example, \citet{KW1}; \citet{Linfinity};
\citet{KulikovLopuhaaJTP}; \citet{PietHL99}; \citet
{DLarxiv}. 
In the misspecified region, the required assumptions are \textup{(P)}~and~\textup{(F)}.
Note also that by Remark~\ref{remL2}, the assumption \textup{(G2)} is required
in the result. Additional discussions of these assumptions, including
directions for future research, are provided in the \citet{suppmat}.



To further illustrate these assumptions, as well as Theorem~\ref{teolinear}, we consider the examples (\ref{lineeg1}) and (\ref
{lineeg2}). In example (\ref{lineeg1}), we have that
%
\begin{eqnarray}
\label{lineeg1KL} \hat f_0(x) &=& 1.5\,1_{[0, 0.5]}(x)+0.5\,1_{(0.5,1]}(x).
\end{eqnarray}
The conditions \textup{(S)} and \textup{(P)} are clearly satisfied, as is \textup{(C)} since
$\mathcal{S}_0 = \mathcal{S}_f$. Lastly, \textup{(F)}~holds with $J=2,
\hat q_1=1.5, \hat q_2=0.5, I_1 = (0, 0.5], I_2 = (0.5, 1]$.\break
Applying Theorem~\ref{teolinear} for $g(x)=x$, we find that
$\bar g(x)=0.75\,1_{[0, 0.5]}(x)+\break 0.25\,1_{(0.5,1]}(x)$, and
$I_2\cap\mathcal{W} = \{a_2, b_2\}$ [hence $\operatorname
{gren}({\mathbb U}^{\mathrm{mod}}_2)=0$]. Therefore,
%
\begin{eqnarray}\label{lineeg1meanlimit}
&& \sqrt{n} \bigl(\hat\mu_n(g)-\hat\mu_0(g) \bigr)\nonumber
\\
&&\qquad \Rightarrow \int_0^1 \bar g(x) \,d{\mathbb
U}_{F_0}(x) + \sqrt{\frac{3}{4}} \int_0^1
\frac{u}{2} \operatorname{gren} \bigl({\mathbb U}_1^{\mathrm{mod}}
\bigr) (u)\,du
\\
&&\qquad = -\frac{1}{2}{\mathbb U}_{F_0}(0.5) + \sqrt{
\frac{3}{16}} \int_0^1 u
\operatorname{gren}({\mathbb U}_1) (u)\,du,\nonumber
\end{eqnarray}
where ${\mathbb U}_{F_0}, {\mathbb U}_1$ are independent Brownian
bridges as defined in Theorem~\ref{teolinear}.
Notably, the limit has a non-Gaussian component.

Example (\ref{lineeg2}) can be analyzed similarly. Here,
\begin{eqnarray*}
\hat f_0(x) &=& 12(x-0.5)^2 1_{[0, 0.25]}(x)+0.75\,1_{(0.25,1]}(x).
\end{eqnarray*}
Again, the conditions \textup{(S)} and \textup{(P)} clearly hold. On $S_c = [0, 0.25]$,
we have $\sup_{x \in S_c}|f'_0(x)|=12$ and, therefore, condition \textup{(C)} holds.
On $\mathcal{S}_f = (0.25, 1]$ we have $J=1$, and hence \textup{(F)} also holds.
Applying Theorem~\ref{teolinear} for $g(x)=x$,
we find that $\bar g(x)=x 1_{[0, 0.25]}(x)+ (5/8) 1_{(0.25,1]}(x)$,
and $I_1\cap\mathcal{W} = \{a_1, b_1\}$ [hence $\operatorname
{gren}({\mathbb U}^{\mathrm{mod}}_1)=0$]. Therefore,
\begin{eqnarray*}
\sqrt{n} \bigl(\hat\mu_n(g)-\hat\mu_0(g) \bigr) &
\Rightarrow& \int_0^1 \bar g(x) \,d{\mathbb
U}_{F_0}(x).
\end{eqnarray*}
That is, the limit is zero-mean Gaussian with variance $\sigma^2
\approx0.07032$.



%
\begin{rem}\label{remL2}
Marginal properties of the process $\operatorname{gren}({\mathbb U})$
were studied in
\citet{CarolanDykstraAOS}. The results include marginal densities
and moments,
including 
$\mathrm{E}[(\operatorname{gren}({\mathbb U})(x))^2] =
0.5(x^2/(1-x)+(1-x)^2/x)$.
It follows that
$\mathrm{E}[\int_0^1(\operatorname{gren}({\mathbb U})(x))^2 \,dx]
=\int_0^1(1-x)^2/x \,dx = \infty$,
and hence 
the limiting process
\begin{eqnarray*}
\bigl\langle g, \operatorname{gren}({\mathbb U}) \bigr\rangle&=& \int
_0^1 g(x) \operatorname{gren}({\mathbb U}) (x)
\,dx,
\end{eqnarray*}
exists only for $g \in L_\beta(\mathcal{S}_f)$ for $\beta>2$.
We would therefore not expect convergence of $\hat\mu_n(g)$ for
$g\in L_\beta(\mathcal{S}_f)$ with $\beta\in[1,2]$.
\end{rem}

\section{Beyond linear functionals: A special case}\label{secentropy}

Entropy measures the\break  amount of disorder or uncertainty in a system and
is closely related to the Kullback--Leibler divergence. Let $T(f) =
\int_0^\infty f(x) \log f(x) \,dx$ denote the entropy functional.
A review of testing and other applications of entropy appears, for
example, in \citet{entropy}.

%
\begin{teo}\label{teoentropy}
Suppose that $\hat f_0$ is bounded, the support of $f_0$ is also
bounded, and that
$f_0/\hat f_0 \leq c_0^2 < \infty$. Then
\begin{eqnarray*}
\sqrt{n} \bigl(T(\hat f_n)-T(\hat f_0) \bigr) &
\Rightarrow& \sigma Z,
\end{eqnarray*}
where $Z$ is a standard normal random variable and
\[
\sigma^2 = \operatorname{var}_{f_0} \bigl(\log \bigl(\hat
f_0(X) \bigr) \bigr) = \operatorname{var}_{\hat f_0} \bigl(\log
\bigl(\hat f_0(X) \bigr) \bigr).
\]
\end{teo}

The proof is made up of two key pieces: (1) tight bounds on the
likelihood ratio from Lemma \ref{lemmalogratio} and (2) specialized
equalities which hold for the Grenander estimator.

%
\begin{lem}\label{lemmalogratio}
Suppose that $\hat f_0$ is bounded, the support of $f_0$ is also
bounded, and that $f_0/\hat f_0 \leq c_0^2 < \infty$. Then
\begin{eqnarray*}
\int\log\frac{\hat f_n}{\hat f_0} \,d{\mathbb F}_n &=& O_p
\bigl(n^{-2/3} \bigr).
\end{eqnarray*}
\end{lem}
We note that the conditions we require here are stronger than those of
\citet{patilea}, Corollary 5.6. However, under those conditions
\citet{patilea} establishes convergence rates on $\int\log\frac
{2\hat f_n}{\hat f_n+\hat f_0} \,d{\mathbb F}_n$, which is\vspace*{-1pt}
not sufficient for our purposes. The condition that $f_0/\hat f_0$
is bounded above was also used in the study of misspecification in
\citet{geerbook}, Section~10.4. The condition that the support
of $f_0$ is bounded is the strongest, whereas the condition that
$\hat f_0$ is bounded may be relaxed somewhat. We discuss this
further in \citet{suppmat}.

\begin{pf*}{Proof of Lemma \ref{lemmalogratio}}
We first show that $\int_0^\infty\varphi(\hat f_n) \,d(\widehat
F_n - {\mathbb F}_n)=0$ for any function~$\varphi$. This follows since
$\widehat F_n(x) \geq{\mathbb F}_n(x)$ with equality at finitely many
touch points, and also $\hat f_n$ is constant between all touch
points. Thus, letting $\tau_1, \tau_2, \ldots, \tau_m$ enumerate
the (random) points of touch, we have
\begin{eqnarray*}
\int_0^\infty\varphi(\hat f_n) \,d(
\widehat F_n - {\mathbb F}_n) &=& \sum
_{i=1}^m \varphi \bigl(\hat f_n(
\tau_i) \bigr) \bigl((\widehat F_n - {\mathbb
F}_n) (\tau_i)-(\widehat F_n - {\mathbb
F}_n) (\tau_{i-1}) \bigr) = 0
\end{eqnarray*}
with $\tau_0=0$ and $\tau_m=X_{(n)}$. A similar argument also
establishes that
%
\begin{equation}
\label{line} \int_0^\infty\varphi(\hat
f_0) \,d(\widehat F_0 - F_0) =\int
_{\mathcal{M}} \varphi(\hat f_0) \,d(\widehat
F_0 - F_0)=0.\vadjust{\goodbreak}
\end{equation}
For $\varphi(v)=\log v$, it follows that
\begin{eqnarray*}
\sqrt{n} \bigl(T(\hat f_n)-T(\hat f_0) \bigr) &=&
\sqrt{n} \biggl(\int\log\hat f_n \,d\widehat F_n - \int
\log\hat f_0 \,d \widehat F_0 \biggr)
\nonumber
\\
&=& \sqrt{n} \int\log \biggl(\frac{\hat f_n}{\hat f_0} \biggr) \,d{\mathbb
F}_n + \sqrt{n} \int\log\hat f_0 \,d({\mathbb
F}_n-F_0).
\end{eqnarray*}
The first term is $O_p(n^{-1/6})$ by Lemma \ref{lemmalogratio}.
The second term has a Gaussian limit with variance
$\operatorname{var}_{f_0}(\log\hat f_0(X))$. By (\ref{line})
[with $\varphi(v)=\log^2v,\log v$]
this is equal to $\operatorname{var}_{\hat f_0}(\log\hat f_0(X))$.
\end{pf*}

%
\begin{figure}

\includegraphics{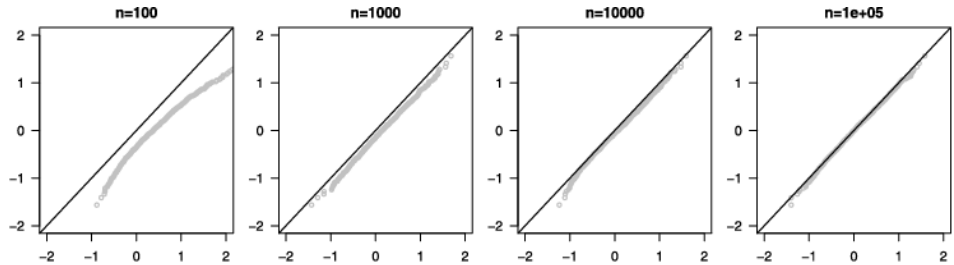}

\caption{Empirical\vspace*{1pt} quantiles of $\sqrt{n}(T(\hat f_n)-T(\hat f_0))$ vs.
the true quantiles of the limiting $N(0,\sigma^2)$ distribution for
$f_0$ given in (\protect\ref{lineeg1}), with $\sigma^2 \approx0.2263$.}
\label{figmissentropy}
\end{figure}

A simulation of this result is shown in Figure~\ref{figmissentropy}
based on the
true density~(\ref{lineeg1}). The KL projection of (\ref{lineeg1})
is given in (\ref{lineeg1KL}). One can easily check that the
conditions of Theorem~\ref{teoentropy} are satisfied in this case.
Note that this density has well-specified flat regions
and, therefore, linear functionals that do not ignore $\mathcal{S}_f
\cap\mathcal{W}$ should
have non-Gaussian terms in their limit; see, for example, (\ref
{lineeg1meanlimit}) for the case when $g(x)=x$. On the other hand, the
entropy functional
will \emph{always} result in a Gaussian limit. The simulations exhibit
a systematic
positive bias. The proof shown above reveals the cause:\vspace*{1pt} the term
$\int\log({\hat f_n}/{\hat f_0} ) \,d{\mathbb F}_n \geq0$
since $\hat f_n$
is the~MLE. In the plots the quantiles of $\sqrt{n}(T(\hat f_n)-T(\hat f_0))$
are shown on the \mbox{$x$-}axis, and these quantiles appear to be shifted to
the right---that is, they are larger than the quantiles of the limiting Gaussian
distribution. 

\section{Conclusion}\label{sec5}

We anticipate that extensions of this work to other one-dimensional
shape-constrained models, such as the log-concave and convex decreasing
constraints, are within reach, although certain technical difficulties
will need to be overcome. In particular, the results of \citet
{patilea} for convex models should yield some results for convex
decreasing densities under misspecification. The Grenander estimator
has a particular simplicity of form, which we have exploited here. Some
progress for the log-concave setting has already been made in
\citet{BJR}, albeit for the discrete (i.e., probability mass
function) setting. We conjecture that statements such as (\ref
{linepower}) will continue to hold for other shape-constraints in
$d=1$ for linear functionals. Similar results for higher-dimensional
shape-constrained models seem premature in view of the current lack of
rate of convergence results even when the model is correctly specified.

\section{Proofs for Section~\texorpdfstring{\protect\ref{seclinear}}{3}}\label{secproofs}

We now present the proof for Theorem \ref{teolinear}.
We proceed by proving convergence results for the
different types of behaviours of the density separately (curved, flat,
misspecified),
and combine the results together at the end. We believe that the intermediate
results are of independent interest to the reader, and we also hope
that this approach
makes the proof more accessible.

\subsection{Strictly curved well-specified density}\label{sec6.1}

We first suppose that the true density $f_0$ satisfies the conditions
introduced in \citet{KW1}.

%
\begin{prop}\label{curved}
Suppose that $f_0$ satisfies conditions \textup{(S)} and \textup{(C)}, and that $g$
satisfies condition \textup{(G1)}. Then
\begin{eqnarray*}
\sqrt{n} \bigl(\hat\mu_n(g)-\mu_0(g) \bigr) &
\Rightarrow& \sigma Z,
\end{eqnarray*}
where $Z$ is a standard normal random variable and $\sigma^2 =
\operatorname{var}(g(X))<\infty$.
\end{prop}

We note that this result is similar to that in \citet{KulikovLopuhaaJTP}.

\begin{pf}
Without loss of generality, we assume that $\mathcal{S}_0=\mathcal
{S}_c=[0,1]$.
Let $\bar{\mu}_n(g)= n^{-1} \sum_{i=1}^n g(X_i)$ denote the empirical
estimator of $\mu_0(g)$. Using Fubini, we have
\begin{eqnarray*}
\bigl|\hat\mu_n(g)-\bar{\mu}_n(g)\bigr| &=& \biggl\llvert\int
_0^1 g'(x) \bigl[\widehat
F_n(x)-{\mathbb F}_n(x) \bigr] \,dx \biggr\rrvert
\\
&\leq& \biggl\{\int
_{\mathcal{S}_c} \bigl|g'(x)\bigr| \,dx \biggr\} \sup
_{x\in\mathcal{S}_c}\bigl| \widehat F_n(x)-{\mathbb
F}_n(x)\bigr|.
\end{eqnarray*}
From the results of \citet{KW1} [see also \citet{DLarxiv},
Corollary 2.2], we have that $\sup_{x\in[0,1]} \sqrt{n}| \widehat
F_n(x)-{\mathbb F}_n(x)| = o_p(n^{-1/6} \log^{2/3} n)$. Therefore,
\begin{eqnarray*}
\sqrt{n} \bigl(\hat\mu_n(g)-\mu_0(g) \bigr) &=&
\sqrt{n} \bigl(\bar\mu_n(g)-\mu_0(g) \bigr)+
o_p \bigl(n^{-1/6} \log^{2/3} n \bigr)
\end{eqnarray*}
from which the result follows.
\end{pf}

\subsection{Piecewise constant well-specified density}\label{sec6.2}

Suppose next that $\mathcal{S}_0 = \mathcal{S}_f = \mathcal{W} \cap
\mathcal{S}_0$.
That is, the true\vadjust{\goodbreak} density is piecewise constant decreasing and can be
expressed as
%
\begin{equation}
\label{defstep} f_0(x) = \sum_{j=1}^J
q_j 1_{(a_j, b_j]} (x),
\end{equation}
where $q_1 > q_2 > \cdots> q_J>0$, $J$ is finite, and $\cup I_j =
\mathcal{S}_0$
where the sets $I_j = (a_j, b_j]$ are disjoint. Indeed, we have
$b_{j}=a_{j+1}$ for $j=1, \ldots, J-1$. 
Note that $p_j = q_j(b_j-a_j)$.
Also, let ${\mathbb U}_1, \ldots, {\mathbb U}_J$ denote independent standard
Brownian bridge processes (each defined on $[0,1]$), and let
$\{Z_1, \ldots, Z_J\}$ be an independent multivariate normal with mean
zero and covariance $\operatorname{diag}(p)-pp^T$ for $p=(p_1, \ldots, p_J)^T$.

%
\begin{prop}\label{step}
Suppose that $f_0$ is as in (\ref{defstep}). Then
$\sqrt{n}(\hat f_n(x)-f_0(x))$ converges weakly to
${\mathbb S}(x)$ in $L_\alpha(\mathcal{S}_f)=L_\alpha(\mathcal
{S}_0)$ for any $\alpha\in[1,2)$, where
\begin{eqnarray*}
{\mathbb S}(x) &=& \sum_{j=1}^J
\frac{1}{b_j-a_j} \biggl\{Z_j + \sqrt{p_j}
\operatorname{gren}({\mathbb U}_j) \biggl(\frac
{x-a_j}{b_j-a_j} \biggr)
\biggr\} 1_{I_j}(x).
\end{eqnarray*}
\end{prop}

A pointwise version of Proposition \ref{step} was originally proved in
\citet{CarolanDykstraCJS}.
Here, we extend these results to convergence in $L_\alpha$, which is a much
stronger statement, requiring tight bounds on the tail behaviour at a
point of the kind
proved in \citet{PietHL99}, Theorem~2.1. In the case of the
decreasing probability
mass function, $\ell_k, k\geq1$ convergence has been established in
\citet{JWejs}.
An immediate corollary of this work is convergence of the linear functionals
$\hat\mu_n(g)$; see Corollary \ref{corstep} below.

\citet{pietuniform}, Theorem 4.1, shows that for $f_0$ equal to
the uniform density on
$[0,1]$ we have
\[
\sqrt{n} \int_0^1 \bigl| \hat{f}_n (x)
-f_0(x) \bigr| \,dx \Rightarrow\int_0^1
\bigl| \operatorname{gren}({\mathbb U}) (x) \bigr| \,dx = 2 \sup_{0 \le x \le1} {
\mathbb U}(x),
\]
where ${\mathbb U}$ is again a standard Brownian bridge process on $[0,1]$.
This is an immediate corollary of Proposition \ref{step} with $J=1$.
On the other hand,
\citet{piet83} [see also \citet{PietPyke}] shows that
\[
\frac{\int_0^1 (\sqrt{n}(\hat f_n(x)-f_0(x)))^2 \,dx - \log
n}{\sqrt{3 \log n}} \Rightarrow Z \sim N(0,1)
\]
and hence convergence of $\sqrt{n}(\hat f_n(x)-f_0(x))$ to
$\operatorname{gren}({\mathbb U})(x)$ in $L_2([0,1])$ fails. See also
Remark~\ref{remL2}.

%
\begin{cor}\label{corstep}
Suppose that $f_0$ takes the form (\ref{defstep}) with bounded support
$\mathcal{S}_0=\mathcal{S}_f\cap\mathcal{W}$ and with $J $
finite.\vadjust{\goodbreak}
Suppose further that
$g$ satisfies condition \textup{(G2)}. Then $\sqrt{n}(\hat\mu_n(g)-\mu
_0(g))\Rightarrow Y_J$, where
\begin{eqnarray*}
Y_J &=& \sum_{j=1}^J \biggl\{
\bar g_j Z_j + \sqrt{p_j}\int
_0^1 g_j(u) \operatorname{gren}({
\mathbb U}_j) (u)\,du \biggr\}
\end{eqnarray*}
with $\bar g_j$ and $g_j$ defined in (\ref{linedefg}).
\end{cor}

In what follows, unless stated otherwise, we assume that $\mathcal{S}_0=[0,1]$.

%
\begin{lem}\label{lemAHstronger}
Suppose that $f_0$ is as in (\ref{defstep}) with a discontinuity at a
point $x_0\neq0$. Then, for all $c>0$,
\[
\sup_{0 \leq x \leq c/n} \bigl\llvert\hat{f}_n(x_0+x)-f_0(x_0+x)
\bigr\rrvert= O_p(1).
\]
\end{lem}

\begin{pf}
It was shown in \citet{AnevskiHossjer}, Theorem~2, that
%
\begin{equation}
\label{lineAHmainresult} \hat{f}_n(x_0+t/n)-\frac{f_0(x_0 -)+f_0(x_0 +)}{2}
\Rightarrow h(t),
\end{equation}
where $h(t)$ is the left derivative of the least concave majorant (over
${\mathbb R}$) of the process
\[
{\mathbb N} \bigl(\lambda(s) \bigr)-\lambda(s)- \biggl\{\frac{f_0(x_0 -)-f_0(x_0
+)}{2}
\biggr\}|s|,
\]
where the rate function is equal to
\begin{eqnarray*}
\lambda(s)&=& \cases{ f_0(x_0 +)s, &\quad$s>0$,
\vspace*{2pt}\cr
f_0(x_0 -)s, &\quad$ s<0$.}
\end{eqnarray*}

Here, ${\mathbb N}$ denotes a standard two-sided Poisson process.
The result in \citet{AnevskiHossjer}, Theorem 2, is established
by a ``switching''
argument similar\vadjust{\goodbreak} to that in the proof of Theorem \ref{teomisslocal}.
The switching argument can also be extended to this situation even if
$f_0(x_0-)= f_0(x_0+)$. A~similar argument may also be used to show
convergence in finite-dimensional distributions as well.
We next show convergence of the supremum norm
\begin{eqnarray*}
&& \sup_{0 \leq x \leq c/n} \bigl\llvert\hat{f}_n
(x_0+x)-f_0(x_0+x) \bigr\rrvert
\\
&&\qquad \Rightarrow \sup_{0 \leq x \leq c} \biggl\llvert h(x)+
\frac{f_0(x_0
-)-f_0(x_0 +)}{2} \biggr\rrvert= O_p(1).
\end{eqnarray*}
This is done by (1) showing that the convergence in (\ref{lineAHmainresult})
also holds in $D[0,\infty)$, and (2) showing that this implies
convergence of the supremum norm (as above).
Both of these steps follow exactly the same argument as the proof of
Theorem 1.1 in
\citet{BJPSW10}, and we therefore omit the details.
\end{pf}

%
\begin{lem}\label{lemexpbounds}
Suppose that $f_0$ is decreasing on ${\mathbb R}$ and flat on $(a,b]$
and fix $x\in(a,b)$.
Then, for any $t_0>0$ and $k_0>0$, there exists a constant
$c_0=t_0/(f_0(b)+t_0/k_0)$ such that
\begin{eqnarray*}
P \bigl(\hat{f}_n(x) > f_0(x)+n^{-1/2}t \bigr)
&\leq& \exp \biggl\{- c_0\frac{t(x-a)}{2} \biggr\} \qquad\mbox{for
all } t\geq t_0
\end{eqnarray*}
for all $n\geq(k_0/3)^2$. Also,
\[
P \bigl(\hat{f}_n(x) < f_0(x)-n^{-1/2}t \bigr)
\leq \exp \biggl\{- \frac
{t^2(b-x)}{2f_0(b)} \biggr\}
\qquad \mbox{for all } t \in \bigl[0, \sqrt{n}f_0(x) \bigr]
\]
and otherwise the probability is equal to zero.
\end{lem}

\begin{pf}
Let ${\mathbb F}_n(a,s)={\mathbb F}_n(s)-{\mathbb F}_n(a)$, and we write
$f_0(a+)= \lim_{x\rightarrow a_+} f_0(x)$. By the switching relation,
\begin{eqnarray*}
&& P \bigl(\hat{f}_n(x) > f_0(x)+n^{-1/2}t
\bigr)
\\
&&\qquad = P \Bigl(\mathop{\operatorname{argmax}}_{s\in[0,1]} \bigl\{{\mathbb
F}_n(s)- \bigl(f_0(x)+n^{-1/2}t \bigr)s \bigr
\}>x \Bigr)
\\
&&\qquad = P \Bigl(\mathop{\operatorname{argmax}}_{s\in[0,1]} \bigl\{{\mathbb
F}_n(a,s)- \bigl(f_0(a+)+n^{-1/2}t \bigr) (s-a)
\bigr\}>x \Bigr)
\\
&&\qquad \leq P \bigl(\bigl\{{\mathbb F}_n(a,s)\geq \bigl(f_0(a+)+n^{-1/2}t
\bigr) (s-a)\bigr\},\mbox{ for some } s\in(x,1] \bigr)
\\
&&\qquad = P \biggl(\frac{{\mathbb F}_n(a,s)}{F_0(a,s)}\geq\frac
{(f_0(a+)+n^{-1/2}t)(s-a)}{F_0(a,s)}, \mbox{ for some } s
\in(x,1] \biggr)
\\
&&\qquad = P \biggl(\frac{{\mathbb F}_n(a,s)}{F_0(a,s)}\geq1+\frac
{n^{-1/2}t}{f_0(a+)}, \mbox{ for some } s\in(x,1] \biggr)
\\
&&\qquad \leq P \biggl(\sup_{s\in(x,1]}\frac{{\mathbb
F}_n(a,s)}{F_0(a,s)}\geq1+
\frac{n^{-1/2}t}{f_0(a+)} \biggr).
\end{eqnarray*}
Since ${\mathbb F}_n(a,s)$ is a binomial random variable, we can bound
the above
using \citet{ShorackWellnerBook}, Inequality 10.3.2, page~416, with
$h(v) = v(\log v-1) +1$ and $\psi(v)= 2 h(1+v)/v^2 \geq(1+v/3)^{-1}$.
It therefore follows that
\begin{eqnarray*}
&& P \biggl(\sup_{s\in(x,1]}\frac{{\mathbb
F}_n(a,s)}{F_0(a,s)}\geq1+
\frac{n^{-1/2}t}{f_0(a+)} \biggr)
\\
&&\qquad \leq \exp \biggl\{-n F_0(a,x) h \biggl(1+
\frac
{n^{-1/2}t}{f_0(a+)} \biggr) \biggr\}
\\
&&\qquad = \exp \biggl\{- \frac{t^2 (x-a)}{2f_0(a+)} \psi \biggl(\frac
{n^{-1/2}t}{f_0(a+)}
\biggr) \biggr\}
\\
&&\qquad \leq \exp \biggl\{- \frac{t(x-a)}{2} \frac
{t/f_0(a+)}{1+(t/f_0(a+))/(3\sqrt{n})} \biggr\}.
\end{eqnarray*}
Write $u=t/f_0(a_+)$ and note that for all $n\geq(k_0/3)^2$ we have
\begin{eqnarray*}
\frac{u}{1+u/(3\sqrt{n})}&\geq& \frac{u}{1+u/k_0},
\end{eqnarray*}
which is a increasing function of $u$. Fix $t_0>0$ and let $u_0=t_0/f_0(a+)$.
Then, with $c_0=u_0/(1+u_0/k_0)=t_0/(f_0(a+)+t_0/k_0)$ we have that
\begin{eqnarray*}
P \bigl(\hat{f}_n(x) > f_0(x)+n^{-1/2}t \bigr)
&\leq& \exp \biggl\{ -c_0 t \biggl(\frac{x-a}{2} \biggr) \biggr
\}.
\end{eqnarray*}

We handle the other side in a similar manner.
\begin{eqnarray*}
&& P \bigl(\hat{f}_n(x) < f_0(x)-n^{-1/2}t
\bigr)
\\
&&\qquad = P \Bigl(\mathop{\operatorname{argmax}}_{s\in[0,1]} \bigl\{{\mathbb
F}_n(s,b)- \bigl(f_0(b)-n^{-1/2}t \bigr) (x-b)
\bigr\}<x \Bigr)
\\
&&\qquad \leq P \biggl(\frac{{\mathbb F}_n(s, b)}{F_0(s, b)}\leq1-\frac
{n^{-1/2}t}{f_0(b)}, \mbox{ for some } s \in[0,x ) \biggr)
\\
&&\qquad \leq P \biggl(\inf_{s\in[0,x)}\frac{{\mathbb F}_n(s, b)}{F_0(s,
b)}\leq1-
\frac{n^{-1/2}t}{f_0(b)} \biggr).
\end{eqnarray*}
We now bound this using the martingale inequality from \citet
{PietHL99}, Lemma 2.3.
\begin{eqnarray*}
&& P \biggl(\inf_{s\in[0,x)}\frac{{\mathbb F}_n(s, b)}{F_0(s, b)} \leq 1-
\frac{n^{-1/2}t}{f_0(b)} \biggr)
\\
&&\qquad \leq \exp \biggl\{-n F(x,b) h \biggl(1-
\frac
{n^{-1/2}t}{f_0(b)} \biggr) \biggr\}
\\
&&\qquad = \exp \biggl\{- \frac{t^2(b-x)}{2f_0(b)} \psi \biggl(-\frac
{n^{-1/2}t}{f_0(b)} \biggr)
\biggr\}.
\end{eqnarray*}
Now, note that since $\hat{f}_n$ is a density, we only consider $t\leq
\sqrt{n}f_0(b)=\sqrt{n}f_0(x)$.
Therefore, we bound only $\psi(-v)$ for $v\in[0,1]$, for which we
have that \mbox{$\psi(-v)\geq1$.} Thus, it follows that
\begin{eqnarray*}
P \bigl(\hat{f}_n(x) < f_0(x)-n^{-1/2}t \bigr)
&\leq&\exp \biggl\{- \frac
{t^2(b-x)}{2f_0(b)} \biggr\}.
\end{eqnarray*}\upqed
\end{pf}

Let $(x)_+=\max(x,0)$ and $(x)_-=\min(x,0)$.

%
\begin{lem}\label{lemsteppointexp}
Suppose that $f_0$ is flat on $(a,b]$ and fix $x\in(a,b)$, and fix
$\alpha>0$.
Then, there exists a constant $C$ such that
\begin{eqnarray*}
E \bigl[ \bigl\llvert\sqrt{n} \bigl(\hat{f}_n(x)-f_0(x)
\bigr)_- \bigr\rrvert^\alpha \bigr] &\leq& C (b-x)^{-\alpha/2},
\\
E \bigl[ \bigl\llvert\sqrt{n} \bigl(\hat{f}_n(x)-f_0(x)
\bigr)_+ \bigr\rrvert^\alpha \bigr] &\leq& C (x-a)^{-\alpha/2}
\end{eqnarray*}
with the second bound valid only for $(x-a)\geq\tilde c_0/n$, for
some $\tilde c_0>0$.
\end{lem}

\begin{pf}
Using the bounds obtained in Lemma \ref{lemexpbounds}, we find that
\begin{eqnarray*}
&& E \bigl[ \bigl\llvert\sqrt{n} \bigl(\hat{f}_n(x)-f_0(x)
\bigr)_{-} \bigr\rrvert^\alpha \bigr]
\\
&&\qquad = \int
_0^\infty\alpha t^{\alpha-1} P \bigl(\sqrt{n}
\bigl(\hat{f}_n(x)-f_0(x) \bigr)_{-} > t
\bigr) \,dt
\\
&&\qquad = \int_0^{n^{1/2}f_0(x)} \alpha t^{\alpha-1} P
\bigl( \hat{f}_n(x) < f_0(x)-n^{-1/2}t \bigr)
\,dt
\\
&&\qquad \leq \int_0^\infty\alpha t^{\alpha-1} \exp
\biggl\{- \frac
{t^2(b-x)}{2f_0(b)} \biggr\}\,dt
\\
&&\qquad = \Gamma(1+\alpha/2) \biggl(\frac{2f_0(b)}{b-x} \biggr)^{\alpha/2}.
\end{eqnarray*}
For the second inequality, we first fix $t_0>0$. We then have
\begin{eqnarray*}
&& E \bigl[ \bigl\llvert\sqrt{n} \bigl(\hat{f}_n(x)-f_0(x)
\bigr)_{+} \bigr\rrvert^\alpha \bigr]
\\
&&\qquad = \int
_0^\infty\alpha t^{\alpha-1} P \bigl(\sqrt{n}
\bigl(\hat{f}_n(x)-f_0(x) \bigr)_{+} > t
\bigr) \,dt
\\
&&\qquad = \int_0^{t_0} \alpha t^{\alpha-1} P
\bigl( \hat{f}_n(x)>f_0(x)+n^{-1/2} t \bigr) \,dt
\\
&&\quad\qquad{}
+ \int_{t_0}^\infty\alpha t^{\alpha-1} P
\bigl(\hat{f}_n (x)>f_0(x)+n^{-1/2} t \bigr) \,dt
\\
&&\qquad \leq t_0^\alpha+ \int_{t_0}^\infty
\alpha t^{\alpha-1} \exp \biggl\{- c_0\frac{t(x-a)}{2} \biggr
\}\,dt
\\
&&\qquad \leq t_0^\alpha+ \Gamma(\alpha+1) \biggl(
\frac
{2}{c_0(x-a)} \biggr)^{\alpha}.
\end{eqnarray*}
Now, recall that $c_0$ takes the form $t_0/(f_0(b)+t_0/k_0)$.
Therefore, we obtain the bounds
\begin{eqnarray*}
\Gamma(\alpha+1) \biggl(\frac{2}{c_0(x-a)} \biggr)^{\alpha} &\leq&
2^\alpha\Gamma(\alpha+1) (x-a)^{-\alpha} \biggl(\frac
{f_0(b)+t_0/k_0}{t_0}
\biggr)^{\alpha}
\\
&\leq& C_\alpha(1+K)^\alpha \biggl(\frac{f_0(b)}{x-a}
\biggr)^\alpha t_0^{-\alpha}
\end{eqnarray*}
as long as $t_0/k_0 \leq K f_0(b)$ for some choice of $K$.
We optimize the entire quantity in $t_0$ to find that
\begin{eqnarray*}
E \bigl[ \bigl\llvert\sqrt{n} \bigl(\hat{f}_n(x)-f_0(x)
\bigr)_{+} \bigr\rrvert^\alpha \bigr] &\leq& A_\alpha
\biggl(\frac{f_0(b)}{x-a} \biggr)^{\alpha/2}
\end{eqnarray*}
for some new constant $A_\alpha$.
Now, in order for this optimized bound to hold, we need $t_0 \leq K
f_0(b) k_0$, and
\begin{eqnarray*}
K^2 f_0(b)^2k_0^2
&\geq& C_\alpha(1+K)^\alpha \biggl(\frac
{f_0(b)}{x-a}
\biggr)^\alpha.
\end{eqnarray*}
The latter translates to $(x-a) \geq\tilde c_0 n^{-1}$ by using
$k_0^2 \le9n$.
\end{pf}

\begin{pf*}{Proof of Proposition \ref{step}}
The outline of the proof is as follows. We first require pointwise convergence,
which follows from \citet{CarolanDykstraCJS}, Theorem 6.4.
One can also easily extend this to convergence in finite-dimensional
distributions.
The particular form of the limit follows from the following
decomposition of a (time-transformed)
Brownian bridge, which is a generalization of \citet
{ShorackWellnerBook}, Exercise 2.2.11, page 32.
Let $F$ denote any distribution function with compact support, which,
without loss of
generality, we assume to be $[0,1]$. Let $0 = a_1 < b_1=a_2 < \cdots<
b_{J-1}=a_J < b_J=1$.
Let ${\mathbb V}, {\mathbb U}_1, \ldots, {\mathbb U}_J$ denote
independent Brownian bridges. Then
%
\begin{eqnarray}\label{lineBBdecomposition}
\qquad {\mathbb V} \bigl(F(t) \bigr) &\equiv& \sum
_{i=1}^J \Biggl\{ \sum_{j=1}^{i-1}
\Delta{\mathbb V} \bigl(F(a_i) \bigr) + \Delta{\mathbb V}
\bigl(F(a_i) \bigr) \frac{F(t)-F(a_i)}{F(b_i)- F(a_i)}
\nonumber\\[-8pt]\\[-8pt]
&&\hspace*{42pt}{} + \sqrt{F(b_i)-
F(a_i)} {\mathbb U}_i \biggl( \frac{F(t)-F(a_i)}{F(b_i)- F(a_i)}
\biggr) \Biggr\} 1_{(a_i, b_i]} (t),\nonumber
\end{eqnarray}
where $\Delta{\mathbb V}(F(a_j))={\mathbb V}(F(b_j))-{\mathbb V}(F(a_j))$.

Recall that the Grenander operator satisfies $\operatorname
{gren}(a+bt+c h(t))=b+c\operatorname{gren}(h(t))$.
Also note that $F_0$ is linear on $(a_i, b_i]$ by assumption.\vspace*{2pt} Therefore,
from \citet{CarolanDykstraCJS}, the limit of $\sqrt{n}(\hat
f_n(x)-f_0(x))$
at a point $x\in I_i=(a_i, b_i]$ can be written as
\begin{eqnarray*}
\operatorname{gren}_{(a_i, b_i]} \bigl({\mathbb V} \bigl(F_0(t)
\bigr) \bigr) &=& \Delta{\mathbb V} \bigl(F(a_i) \bigr)
\frac{1}{b_i-a_i} + \sqrt{p_i} \operatorname{gren} \biggl({\mathbb
U}_i \biggl(\frac
{t-a_i}{b_i-a_i} \biggr) \biggr)
\\
&=&\frac{1}{b_i-a_i} \biggl\{ \Delta{\mathbb V} \bigl(F(a_i)
\bigr) + \sqrt{p_i} \operatorname{gren} ({\mathbb U}_i )
\biggl(\frac
{t-a_i}{b_i-a_i} \biggr) \biggr\}
\end{eqnarray*}
from the above characterization. Finally,
$\{\Delta{\mathbb V}(F(a_1)), \ldots, \Delta{\mathbb V}(F(a_J))\}
\stackrel{d}{=} \{Z_1, \ldots, Z_J\}$ as in Proposition \ref{step}.

The second step is to show that the process ${\mathbb S}_n(x) = \sqrt {n}(\hat f_n(x)-f_0(x))$
is tight in $L_\alpha(\mathcal{S})$. For this, we first need a
characterization of
compact sets in $L_\alpha(\mathcal{S})$ for $\alpha\geq1$.
These appear, for example, in \citet{Dunford}, page 298 [see also
\citet{simon}].
For $\mathcal{S}$ bounded, a set $\mathcal{K} \subset L_\alpha
(\mathcal{S})$ is relatively compact if for all $f\in\mathcal{K}$:
\begin{longlist}[(2)]
\item[(1)] $\sup_{f\in\mathcal{K}} \int_{\mathcal{S}} |f(x)|^\alpha
\,dx<\infty$,
\item[(2)] $\lim_{\delta\rightarrow0} \sup_{f\in\mathcal{K}} \int_{\mathcal
{S}}|f(x+\delta)-f(x)|^\alpha \,dx \rightarrow0$.
\end{longlist}

We want to show that for each $\epsilon>0$ we can find a compact
subset $\mathcal{K}=\mathcal{K}_{\epsilon}$ of $L_\alpha(\mathcal{S})$
such that
$\limsup_n P( {\mathbb S}_n \in\mathcal{K}^c ) < \epsilon$.
Thus, we want to show that
%
\begin{eqnarray}\label{TightnessPart1}
\limsup_n P \biggl( \int_0^1
\bigl| {\mathbb S}_n (x) \bigr|^{\alpha} \,dx > M \biggr) &\rightarrow&0
\qquad \mbox{as } M\rightarrow\infty
\end{eqnarray}
and
\begin{eqnarray}\label{TightnessPart2}
\limsup_n P \biggl( \int_0^1
\bigl| {\mathbb S}_n (x+\delta) - {\mathbb S}_n (x)
\bigr|^{\alpha} \,dx > \epsilon \biggr) &\rightarrow& 0
\end{eqnarray}
as $\delta\rightarrow0$, for every $\epsilon> 0$.

To show the first of these, we proceed as follows: for $f_0$ as in
(\ref{defstep}),
\[
\int_0^1 \bigl|{\mathbb S}_n (x)
\bigr|^{\alpha} \,dx = \sum_{j=1}^J \int
_{(a_j, b_j]} \bigl| {\mathbb S}_n (x) \bigr|^{\alpha} \,dx
\]
and hence we have
%
\begin{equation}
\label{linetight1} P \biggl( \int_0^1 \bigl| {\mathbb
S}_n (x) \bigr|^{\alpha} \,dx > M \biggr) 
\leq \sum_{j=1}^J P \biggl(
\int_{(a_j, b_j]} \bigl| {\mathbb S}_n (x) \bigr|^{\alpha}
\,dx > M/J \biggr).
\end{equation}
Thus, it suffices to show that
\[
\limsup_n P \biggl( \int_{(a,b]} \bigl| {
\mathbb S}_n (x) \bigr|^{\alpha} \,dx > M \biggr) \rightarrow0
\]
as $M\rightarrow\infty$ for each fixed $(a,b]$ with $f_0$ flat on $(a,b]$.
Now,
%
\begin{eqnarray}\label{DeCompOfIntegralTwoPieces}
P \biggl( \int_{(a,b]} \bigl| {\mathbb S}_n (x)
\bigr|^{\alpha} \,dx > M \biggr) &\le& P \biggl( \int_{(a,a+\tilde{c}_0/n]} \bigl|
{\mathbb S}_n (x) \bigr|^{\alpha} \,dx > M/2 \biggr)
\nonumber\\[-8pt]\\[-8pt]
&& {}+ P \biggl( \int_{(a+ \tilde{c}_0/n,b]} \bigl| {\mathbb S}_n (x)
\bigr|^{\alpha} \,dx > M/2 \biggr) \nonumber
\end{eqnarray}
and we handle each term separately. From Lemma~\ref{lemAHstronger},
it follows that
%
\begin{eqnarray}
\label{lineboundfirstterm} \int_{(a,a+\tilde{c}_0/n]} \bigl| {\mathbb S}_n (x)
\bigr|^{\alpha} \,dx & = & \int_{(a,a+\tilde{c}_0/n]} n^{\alpha/2} \bigl|
\hat{f}_n (x) - f_0 (x) \bigr|^{\alpha} \,dx
\nonumber\\[-8pt]\\[-8pt]
& = & n^{\alpha/2} \tilde{c}_0 n^{-1}
O_p (1) = o_p (1)\nonumber
\end{eqnarray}
for $\alpha< 2$. 
For the second term, we use Markov's inequality,
Lemma~\ref{lemsteppointexp} and Fubini's theorem to get
%
\begin{eqnarray}
\label{lineboundsecondterm} \qquad&& P \biggl( \int_{(a+ \tilde{c}_0/n,b]} \bigl| {\mathbb
S}_n (x) \bigr|^{\alpha} \,dx > M/2 \biggr)
\nonumber
\\
&&\qquad \leq \frac{2}{M} 2^{\alpha-1} \biggl
\{\int_{(a+
\tilde{c}_0/n,b]} E \bigl| {\mathbb S}_n (x)_+
\bigr|^{\alpha} \,dx +\int_{(a+
\tilde{c}_0/n,b]} E \bigl| {\mathbb
S}_n (x)_- \bigr|^{\alpha} \,dx \biggr\}
\\
&&\qquad \le \frac{2^\alpha C}{M} \biggl\{\int_{(a,b]}
(x-a)^{-\alpha/2} \,dx +\int_{(a,b]} (b-x)^{-\alpha/2}
\,dx \biggr\} \leq\widetilde C/M\nonumber 
\end{eqnarray}
for some new, finite, constant $\widetilde C$ depending on $a,b, \alpha
$, noting that $\alpha<2$.
Combining (\ref{lineboundfirstterm}) and (\ref
{lineboundsecondterm}) yields (\ref{TightnessPart1}) for our choice
of $f_0$.

Now, to prove (\ref{TightnessPart2}). Since $f_0 (x)$ is constant for
$x \in(a_j, b_j]$ for each $j$, the processes
${\mathbb S}_n (x) = \sqrt{n} ( \hat{f}_n (x) - f_0 (x))$, are
piecewise monotone,
and hence the convergence in $L_{\alpha} ((a_j, b_j])$ for $ \alpha
\in[1,2)$
and each $j \le J$ follows as
in \citet{MR1311975}, Corollary 2, page 1260. We conclude that
(\ref{TightnessPart2}) holds, and hence ${\mathbb S}_n$ is tight in
$L_{\alpha} ({\mathcal S})$ when $\alpha<2$.
\end{pf*}



\begin{pf*}{Proof of Corollary \ref{corstep}}
Convergence follows immediately by continuity of the linear functional
$\int g(x) {\mathbb S}_n(x)\,dx$ by H\"{o}lder's inequality.
We need only check the final form, that is, $\int g(x) {\mathbb
S}(x)\,dx$ is equal to
\begin{eqnarray*}
\hspace*{-10pt}&& \sum_{i=1}^J \biggl
\{Z_i \frac{\int_{a_i}^{b_i}
g(x)\,dx}{b_i-a_i} +\sqrt{p_i}\int
_{a_i}^{b_i} \frac{g(x)}{b_i-a_i}\operatorname{gren}({
\mathbb U}_i) \biggl(\frac{x-a_i}{b_i-a_i} \biggr)\,dx \biggr\}
\\
\hspace*{-10pt}&&\qquad = \sum_{j=1}^J \biggl
\{Z_i \frac{\int_{a_i}^{b_i} g(x)\,dx}{b_i-a_i} +\sqrt{p_i}\int
_0^1 g \bigl((b_i-a_i)u+a_i
\bigr) \operatorname{gren}({\mathbb U}_j) (u )\,du \biggr\}.
\end{eqnarray*}\upqed
\end{pf*}

\subsection{Piecewise constant KL density}\label{sec6.3}

We next consider the case that $\hat f_0(x)$ can be written in
the form (\ref{defKLform}) with condition \textup{(F)}. Let ${\mathbb U}_1,
\ldots, {\mathbb U}_J$, denote independent
standard Brownian bridge processes (each defined on $[0,1]$), and for each
$j$ define ${\mathbb U}_j^{\mathrm{mod}}$ as in Remark \ref{remmisslocal} with
$I_j = [a_j, b_j]$ replacing $[a,b]$. Also, let $\{Z_1, \ldots, Z_J\}$
be an independent
multivariate normal with mean zero and covariance
$\operatorname{diag}(p)-pp^T$ for $p=(p_1, \ldots, p_J)^T$, where
$p_j = F_0(b_j)-F_0(a_j)=\widehat F_0(b_j)-\widehat F_0(a_j)$.

%
\begin{prop}\label{stepmiss}
Suppose that $\hat f_0$ satisfies conditions \textup{(P)} and \textup{(F)} with
$\mathcal{S}_0=\mathcal{S}_f$ and that $g$ satisfies condition \textup{(G2)}.
Then $\sqrt{n}(\hat f_n(x)-f_0(x))$ converges weakly  to
${\mathbb S}^{\mathrm{mod}}(x)$ in $L_\alpha(\mathcal{S})$ for $\alpha\in
(0,2)$, where
\begin{eqnarray*}
{\mathbb S}^{\mathrm{mod}}(x) &=& \sum_{j=1}^J
\frac{1}{b_j-a_j} \biggl\{Z_j + \sqrt{p_j}
\operatorname{gren} \bigl({\mathbb U}_j^{\mathrm{mod}} \bigr) \biggl(
\frac
{x-a_j}{b_j-a_j} \biggr) \biggr\} 1_{(a_j, b_j]}(x).
\end{eqnarray*}
\end{prop}

%
\begin{cor}\label{corstepmiss}
Suppose that $\hat f_0$ satisfies conditions \textup{(P)} and \textup{(F)}
with \mbox{$\mathcal{S}_0=\mathcal{S}_f$} and that $g$ satisfies condition \textup{(G2)}.
Then $\sqrt{n}(\hat\mu_n(g)-\hat\mu_0(g))\Rightarrow Y_J$, where
\begin{eqnarray*}
Y_J^{\mathrm{mod}} &=& \sum_{j=1}^J
\biggl\{ \bar g_j Z_j + \sqrt{p_j}\int
_0^1 g_j(u) \operatorname{gren}
\bigl({\mathbb U}_j^{\mathrm{mod}} \bigr) (u)\,du \biggr\}
\end{eqnarray*}
with $g_j(u)$ and $\bar g_j$ defined as in (\ref{linedefg}).
\end{cor}

The proof of these results is very close to that of Proposition \ref{step},
and we omit any details which are the same. The difference lies in the
following modifications to Lemmas \ref{lemAHstronger} and \ref
{lemexpbounds}.
Note that we add the additional requirement that $f_0$ be bounded below \textup{(P)}.

%
\begin{lem}\label{lemAHstrongermiss}
Fix a point $x\in\mathcal{S}_0$ and let $[a,b]$ denote the largest
interval $I$
such that $x\in I$ and $\hat f_0$ is constant on $I$. Then, for all $c>0$,
\[
\sup_{0 \leq u \leq c/n} \bigl\llvert\hat{f}_n(a+u)-\hat
f_0(a+u) \bigr\rrvert= O_p(1).
\]
\end{lem}

\begin{pf}
By the switching relation, it follows that
\begin{eqnarray*}
&&P \bigl(\hat{f}_n(a+u/n)-\hat f_0(a+u/n)<t \bigr)
\\
&&\qquad = P \bigl(\hat{f}_n(a+u/n)-\hat f_0(b)<t
\bigr)
\\
&&\qquad =P \Bigl(\mathop{\operatorname{argmax}}_{z\in[0,1]} \bigl\{{\mathbb
F}_n(z)- \bigl(\hat f_0(b)+t \bigr)z \bigr\}<a+u/n
\Bigr)
\\
&&\qquad =P \Bigl(n \Bigl(\mathop{\operatorname{argmax}}_{z\in[0,1]} \bigl
\{{\mathbb F}_n(z)- \bigl(\hat f_0(b)+t \bigr)z \bigr\}-a
\Bigr)<u \Bigr)
\end{eqnarray*}
and the inner process
\begin{eqnarray*}
&&n \Bigl(\mathop{\operatorname{argmax}}_{z\in[0,1]} \bigl\{ {\mathbb
F}_n(z)- \bigl(\hat f_0(b)+t \bigr)z \bigr\}-a \Bigr)
\\
&&\qquad =\mathop{\operatorname{argmax}}_{h\geq-n a} \bigl\{{\mathbb
F}_n(a+h/n)- \bigl(t+\hat f_0(b) \bigr) (a+h/n) \bigr\}
\\
&&\qquad = \mathop{\operatorname{argmax}}_{h\geq-n a} \bigl\{{\mathbb
V}_n(h) \bigr\}, 
\end{eqnarray*}
where
${\mathbb V}_n(h)=A_n(h)+B_n(h)-th$, with $A_n(h)
=n({\mathbb F}_n(a+h/n)-{\mathbb F}_n(a))-n(F_0(a+h/n)-F_0(a))$
and $B_n(h)=n(F_0(a+h/n)-F_0(a)) - \hat f_0(b)h$.

Now, the term ${\mathbb N}_n(h)=n({\mathbb F}_n(a+h/n)-{\mathbb
F}_n(a))$ is binomial with
mean $n(F_0(a+h/n)-F_0(a)) \rightarrow f_0(a) h$. Therefore, $A_n(h)$
converges to a
centered\vadjust{\goodbreak} Poisson random variable with mean $f_0(a)$.
A similar argument may be used to show convergence as a process of
$A_n(h)\Rightarrow{\mathbb N}(h)-f_0(a)h$,
where ${\mathbb N}(\cdot)$ is a Poisson process with rate $\lambda(h)=f_0(a)$.
The second piece, $B_n(h)$ satisfies
\begin{eqnarray*}
n^{-1}B_n(h) \cases{ = 0, &\quad$h \in n \bigl\{[a,b]\cap
\mathcal{W} - a \bigr\}$,
\vspace*{2pt}\cr
< 0, &\quad$h \in n \bigl\{[a,b]\cap\mathcal{M} - a
\bigr\}$.}
\end{eqnarray*}
Thus, if for all $\delta>0$ $[a,\delta) \cap\mathcal{W} = \{a\}$
then the limit of
$B_n(h)$ is $0$ if $h=0$ and is equal to $-\infty$ otherwise [we will call
this setting case (A)]. If the above assumption is not true [we will
call this setting case (B)],
then $\lim_{n\rightarrow\infty}B_n(h)=0$ for all $h\geq0$. In case
(A), it follows that the
limit of ${\mathbb V}_n(h)$ is equal to 0 at $h=0$ and is equal to
$-\infty$ otherwise.
Therefore, $\operatorname{argmax}_{h\geq0}\{{\mathbb V}_n(h)\}=0$
here. In case (B), the limit of
${\mathbb V}_n(h)$ is a centered (a.k.a. compensated) Poisson process
with rate
$f_0(a)$. We therefore have that, in case (A),
\begin{eqnarray*}
&& P \bigl(\hat{f}_n(a+u/n)-\hat f_0(a+u/n)<t \bigr)
\\
&&\qquad
= P \Bigl(\mathop{\operatorname{argmax}}_{h\geq{-na}} \bigl\{{
\mathbb V}_n(h) \bigr\}<u \Bigr) \rightarrow1
\end{eqnarray*}
and in case (B),
\begin{eqnarray*}
&& P \bigl(\hat{f}_n(a+u/n)-\hat f_0(a+u/n)<t \bigr)
\\
&&\qquad =
P \Bigl(\mathop{\operatorname{argmax}}_{h\geq{-na}} \bigl\{{\mathbb
V}_n(h) \bigr\}<u \Bigr)
\\
&&\qquad \rightarrow P \Bigl(\mathop{\operatorname{argmax}}_{h\geq0}
\bigl\{{\mathbb N}(h)-f_0(a)h-th \bigr\}<u \Bigr),
\end{eqnarray*}
which\vspace*{1pt} gives us pointwise convergence in distribution in both cases.

Lastly, note that $\hat f_0(a+u)=\hat f_0(b)$ is a constant, and
$\hat f_n(a+u)$ is decreasing in~$u$ by definition. Therefore,
$\sup_{0\leq u\leq c/n}|\hat f_n(a+u)-\hat f_0(a+u)| =
|\hat f_n(a)-\hat f_0(b)|$,
which converges as described above.
\end{pf}

%
\begin{lem}\label{lemexpboundsmiss}
Suppose that $\hat f_0$ is flat on $(a,b]$ and fix $x\in(a,b)$.
Assume also that $\inf_{x\in(a,b]} f_0(x) = \alpha_0>0$, and let
$\hat c_0 = \alpha_0/\hat f_0(b)$. Then, for any $t_0>0$ and $k_0>0$,
there exists a constant $c_0=t_0/(\hat f_0(b)+t_0/k_0)$ such that
\begin{eqnarray*}
P \bigl(\hat{f}_n(x) > \hat f_0(x)+n^{-1/2}t
\bigr) &\leq& \exp \biggl\{- \hat c_0 c_0
\frac{t(x-x_0)}{2} \biggr\} \qquad\mbox{for all } t\geq t_0
\end{eqnarray*}
for all $n\geq(k_0/3)^2$. Also, for all $t \in[0, \sqrt{n}\hat f_0(x)]$,
\begin{eqnarray*}
P \bigl(\hat{f}_n(x) < \hat f_0(x)-n^{-1/2}t
\bigr) &\leq& \exp \biggl\{- \hat c_0 \frac{t^2(b-x)}{2\hat f_0(x)} \biggr\}
\end{eqnarray*}
and otherwise the probability is equal to zero.
\end{lem}

\begin{pf}
Let ${\mathbb F}_n(a,s)={\mathbb F}_n(s)-{\mathbb F}_n(a)$, and we write
$\hat\theta= \hat f_0(x)$.
Repeating the argument for the proof of Lemma \ref{lemexpbounds}, we
obtain that
\begin{eqnarray*}
&& P \bigl(\hat{f}_n(x) > \hat f_0(x)+n^{-1/2}t
\bigr)
\\[-2pt]
&&\qquad \leq P \biggl(\frac{{\mathbb F}_n(a,s)}{F_0(a,s)}\geq\frac
{(\hat\theta+n^{-1/2}t)(s-a)}{F_0(a,s)}, \mbox{ for some }
s\in(x,1] \biggr)
\\[-2pt]
&&\qquad \leq P \biggl(\sup_{s\in(x,1]}\frac{{\mathbb F}_n(a,s)}{F_0(a,s)} \geq \inf
_{s\in(x,1]}\frac{(\hat\theta
+n^{-1/2}t)(s-a)}{F_0(a,s)} \biggr)
\\[-2pt]
&&\qquad \leq P \biggl(\sup_{s\in(x,1]}\frac{{\mathbb
F}_n(a,s)}{F_0(a,s)}\geq\inf
_{s\in(x,1]}\frac{(\hat\theta
+n^{-1/2}t)(s-a)}{\widehat F_0(a,s)} \biggr)
\\[-2pt]
&&\qquad = P \biggl(\sup_{s\in(x,1]}\frac{{\mathbb F}_n(a,s)}{F_0(a,s)}\geq1+
\frac{n^{-1/2}t}{\hat\theta} \biggr)
\end{eqnarray*}
since $\widehat F_0(s) > F_0(s)$ with equality at $s=a,b$.
Applying the exponential bounds for binomial variables as before, we
find that
\begin{eqnarray*}
&& P \biggl(\sup_{s\in(x,1]}\frac{{\mathbb
F}_n(a,s)}{F_0(a,s)}\geq1+
\frac{n^{-1/2}t}{\hat\theta} \biggr)
\\[-2pt]
&&\qquad \leq \exp \biggl\{-n F_0(a,x) h \biggl(1+
\frac{n^{-1/2}t}{\hat
\theta} \biggr) \biggr\}
\\[-2pt]
&&\qquad \leq \exp \biggl\{- \biggl[\inf
_{x\in[a,b]} \frac{f_0(x)}{\hat
\theta} \biggr] \frac{t(x-a)}{2}
\frac{t/\hat\theta}{1+(t/\hat\theta
)/(3\sqrt{n})} \biggr\}.
\end{eqnarray*}
Therefore, assuming that $\inf_{x\in[a,b]}\frac{f_0(x)}{\hat\theta
} = \hat c_0 >0$,
we can repeat the same argument as for Lemma (\ref{lemexpbounds}).

We handle the other side in a similar manner:
\begin{eqnarray*}
P \bigl(\hat{f}_n(x) < \hat f_0(x)-n^{-1/2}t
\bigr)
&\leq& P \biggl(\inf_{s\in[0,x)}\frac{{\mathbb
F}_n(s,b)}{F_0(s,b)}\leq\sup
_{s\in[0,x)}\frac{(\hat\theta
-n^{-1/2}t)(s-b)}{F_0(s,b)} \biggr)
\\[-2pt]
&\leq &P \biggl(\inf_{s\in[0,x)}\frac{{\mathbb
F}_n(s,b)}{F_0(s,b)}\leq1-
\frac{n^{-1/2}t}{\hat\theta} \biggr).
\end{eqnarray*}
We again bound this using the martingale inequality from \citet
{PietHL99}, Lemma 2.3:
\begin{eqnarray*}
&& P \biggl(\inf_{s\in[0,x)}\frac{{\mathbb F}_n(s,
b)}{F_0(s, b)} \leq1-
\frac{n^{-1/2}t}{\hat\theta} \biggr)
\\[-2pt]
&&\qquad \leq \exp \biggl\{-n F(x,b) h \biggl(1-\frac{n^{-1/2}t}{\hat\theta
} \biggr)
\biggr\}
\\[-2pt]
&&\qquad = \exp \biggl\{- \biggl[ \inf_{x\in[a,b]}
\frac{f_0(x)}{\hat\theta
} \biggr] \frac{t^2(b-x)}{2\hat\theta} \psi \biggl(-\frac
{n^{-1/2}t}{\hat\theta}
\biggr) \biggr\}.
\end{eqnarray*}\upqed
\end{pf}

\subsection{Putting it all together}\label{sec6.4}
\mbox{}

\begin{pf*}{Proof of Theorem \ref{teolinear}}
To illustrate the method of proof, we consider a simplified case.
Since $g \in L_{\beta} ({\mathcal S}_f)$ for some $\beta>2$
and ${\mathbb S}_n$ converges in $L_{\alpha} ({\mathcal S})$ the proof easily
extends to a general setting. Suppose then that $\mathcal{S}_c=[0,a]$
and $\mathcal{S}_f=[a,b]$, so that the support is $\mathcal{S}_0 = [0,b]$.
Furthermore, we assume that on $\mathcal{S}_f$ we have $J=1$.
Let ${\mathbb U}_n(x)=\sqrt{n}({\mathbb F}_n(x)-F_0(x))$. Then
\begin{eqnarray*}
&& \sqrt{n} \bigl(\hat\mu_n(g)-\hat\mu_0(g) \bigr)
\\
&&\qquad = \int_0^a g(x) \,d{\mathbb
U}_n(x) + \int_a^b g(x) \sqrt{n}
\bigl(\operatorname{gren}({\mathbb F}_n) (x)-\hat f_0(x)
\bigr) \,dx + \varepsilon_n,
\end{eqnarray*}
where $\varepsilon_n = \sqrt{n}\int_0^a g(x) \,d(\widehat F_n-{\mathbb
F}_n)(x)$.
From assumptions \textup{(C)} and \textup{(G1)}, it follows that $\varepsilon_n =
o_p(1)$ as in Proposition \ref{curved}.

Next, let ${\mathbb W}_n= {\mathbb U}_n(b)-{\mathbb U}_n(a)$, and let
${\mathbb V}_n(x)= {\mathbb U}_n(x)-{\mathbb U}_n(a)-\frac
{x-a}{b-a}{\mathbb W}_n$ for $x\in[a,b]$.
Lastly, let $\ell(x)=F_0(a)+\hat f_0(b)(x-a)$. Then for $x\in(a,b]$,
\begin{eqnarray*}
&& \sqrt{n} \bigl(\operatorname{gren}({\mathbb F}_n) (x) - \hat
f_0(x) \bigr)
\\
&&\qquad = \sqrt{n} \bigl(\operatorname{gren}({\mathbb
F}_n) (x) - \hat f_0(b) \bigr)
=\operatorname{gren} \bigl({\mathbb U}_n + \sqrt{n}(F_0-
\ell) \bigr)
\\
&&\qquad
= \frac{1}{b-a} {\mathbb W}_n + \operatorname{gren} \bigl({
\mathbb V}_n + \sqrt{n}(F_0-\ell) \bigr)
\end{eqnarray*}
and we also define ${\mathbb V}_n^{\mathrm{mod}}={\mathbb V}_n + \sqrt {n}(F_0-\ell)$. Therefore, $\sqrt{n}(\hat\mu_n(g)-\hat\mu
_0(g))$ is equal to 
\begin{eqnarray*}
&& \int_0^a g(x)\,d{ \mathbb
U}_n(x) + {\mathbb W}_n \frac
{1}{b-a} \int
_a^b g(x) + \int_a^b
g(x) \operatorname{gren} \bigl({\mathbb V}_n^{\mathrm{mod}} \bigr)
(x)\,dx +o_p(1)
\\
&&\qquad = \int_0^b \bar g(x) \,d{\mathbb
U}_n(x)+ \int_a^b g(x)
\operatorname{gren} \bigl({\mathbb V}_n^{\mathrm{mod}} \bigr) (x)\,dx
+o_p(1)
\end{eqnarray*}
from the definition of ${\mathbb W}_n$ and of $\bar g$. The weak
limit of ${\mathbb V}_n^{\mathrm{mod}}$
can be established similarly as in Theorem \ref{teomisslocal} and
Remark \ref{remmisslocal}.
The outline of the rest of the proof proceeds as follows:
\begin{longlist}[(3)]
\item[(1)]
Joint weak convergence of
$\{\int_0^b \bar g\, d{\mathbb U}_n, {\mathbb V}_n^{\mathrm{mod}}(x_1), \ldots,
{\mathbb V}_n^{\mathrm{mod}}(x_k)\}$ to a Gaussian limit.

\item[(2)]
Joint weak convergence of
\[
\biggl\{\int_0^b \bar g\, d{\mathbb U}_n, \operatorname{gren}\bigl({\mathbb
V}_n^{\mathrm{mod}}\bigr)(x_1), \ldots, \operatorname{gren}\bigl({\mathbb
V}_n^{\mathrm{mod}}\bigr)(x_k)\biggr\}
\]
via the switching relation.

\item[(3)]
We have that
\begin{eqnarray*}
\operatorname{gren} \bigl({\mathbb V}_n^{\mathrm{mod}} \bigr) (x)
&=& \sqrt{n} \bigl(\hat f_n(x)-\hat f_0(x) \bigr) -
\frac
{1}{b-a}{\mathbb W}_n,
\end{eqnarray*}
where in Proposition \ref{stepmiss} we showed that the first term on
the right-hand side is tight in $L_\alpha(a,b)$. The second term on the
right-hand side
is a tight constant and, therefore, $\operatorname{gren}({\mathbb
V}_n^{\mathrm{mod}})(x)$ is also tight in $L_\alpha(a,b)$.

\item[(4)] From (1) and (3), we obtain marginal tightness of the terms
$\int_0^b \bar g\, d{\mathbb U}_n$ in ${\mathbb R}$ and $\operatorname
{gren}({\mathbb V}_n^{\mathrm{mod}})(\cdot)$ in
$L_\alpha(a,b)$, which implies joint tightness in ${\mathbb R}\times
L_\alpha$.
The full result now follows by the continuous mapping theorem.
\end{longlist}

Lastly, we note that since $F_0(z)=\widehat F_0(z)$ at $z=a,b$ and
$\bar g$ is constant on $[a,b]$ then
$\int_{\mathcal{S}_0} \bar g(x) \,d{\mathbb U}_{F_0}(x)=\int_{\mathcal
{S}_0} \bar g(x) \,d{\mathbb U}_{\widehat F_0}(x)$.
\end{pf*}

\section*{Acknowledgements}
The author thanks Valentin Patilea for sharing a copy of his thesis,
Takumi Saegusa for pointing out a small error in one of the proofs and
the referees for a number of helpful suggestions.
Parts of this work were completed while the author was visiting the
University of
Washington and the University of Heidelberg, and the author thanks both
institutions for their hospitality and financial travel support, and in
particular, Tilmann Gneiting for MATCH funding. The author also thanks
Jon Wellner for generous contributions to this work.

\begin{supplement}
\stitle{Supplement to ``Convergence of linear functionals of the Grenander estimator under misspecification''}
\slink[doi]{10.1214/13-AOS1196SUPP} 
\sdatatype{.pdf}
\sfilename{AOS1196\_supp.pdf}
\sdescription{We provide some proofs and technical details, as well as
additional discussions of the assumptions in Theorems~\ref{teolinear}
and \ref{teoentropy}.}
\end{supplement}



%

\printaddresses

\end{document}